%% file: main.tex
\documentclass[12pt,a4paper,leqno]{article}

\usepackage{latexsym}
\usepackage{amssymb,exscale}
\usepackage[centertags]{amsmath}
\usepackage{amsthm}

\usepackage[dvips]{hyperref}

\input{preamble}

\begin{document}

\title{Moderate deviations for random fields and random complex zeroes}

\author{Boris Tsirelson}

\date{}
\maketitle

\stepcounter{footnote}
\footnotetext{%
 This research was supported by \textsc{the israel science foundation}
 (grant No.~683/05).}

\begin{abstract}
Moderate deviations for random complex zeroes are deduced from a new theorem
on moderate deviations for random fields.
\end{abstract}

\setcounter{tocdepth}{1}
\tableofcontents

\section*{Introduction}
\input{intro}

\numberwithin{equation}{section}

\section[A chain of inequalities]
  {\raggedright A chain of inequalities}
\label{sect1}
\input{sect1}

\section[Splittable random processes]
  {\raggedright Splittable random processes}
\label{sect2}
\input{sect2}

\section[Logarithm of a Gaussian process]
  {\raggedright Logarithm of a Gaussian process}
\label{sect3}
\input{sect3}

\section[An example of a splittable process]
  {\raggedright An example of a splittable process}
\label{sect4}
\input{sect4}

\section[Splittable random fields]
  {\raggedright Splittable random fields}
\label{sect5}
\input{sect5}

\section[An example of a splittable field]
  {\raggedright An example of a splittable field}
\label{sect6}
\input{sect6}

\section[Random complex zeroes]
  {\raggedright Random complex zeroes}
\label{sect7}
\input{sect7}

\bigskip
\filbreak
{
\small
\begin{sc}
\parindent=0pt\baselineskip=12pt
\parbox{4in}{
Boris Tsirelson\\
School of Mathematics\\
Tel Aviv University\\
Tel Aviv 69978, Israel
\smallskip
\par\quad\href{mailto:tsirel@post.tau.ac.il}{\tt
 mailto:tsirel@post.tau.ac.il}
\par\quad\href{http://www.tau.ac.il/~tsirel/}{\tt
 http://www.tau.ac.il/\textasciitilde tsirel/}
}

\end{sc}
}
\filbreak

\end{document}

%% file: preamble.tex
\swapnumbers
\theoremstyle{definition}
\newtheorem{definition}[equation]{Definition}
\newtheorem{theorem}[equation]{Theorem}
\newtheorem{proposition}[equation]{Proposition}
\newtheorem{lemma}[equation]{Lemma}
\newtheorem{corollary}[equation]{Corollary}
\newtheorem{remark}[equation]{Remark}

\renewcommand{\phi}{\varphi}

\newcommand{\I}{{\rm i}}
\newcommand{\D}{\mathrm{d}}
\newcommand{\E}{\mathrm{e}}

\renewcommand{\(}{\bigl(}
\renewcommand{\)}{\bigr)\vphantom{)}}

\renewcommand{\Re}{\operatorname{Re}}
\renewcommand{\Im}{\operatorname{Im}}

\newcommand{\ip}[2]{\langle#1,#2\rangle}
\newcommand{\imply}{\;\;\;\Longrightarrow\;\;\;}
\newcommand{\imp}{$ \Longrightarrow $ }
\newcommand{\pd}{\partial}

\newcommand{\ti}{\tilde}

\newcommand{\sgn}{\operatorname{sgn}}
\newcommand{\const}{\operatorname{const}}

\newcommand{\One}{\mathbf1}

\newcommand{\old}{\text{old}}
\newcommand{\new}{\text{new}}

\newcommand{\eps}{\varepsilon}
\newcommand{\si}{\sigma}
\newcommand{\ga}{\gamma}
\newcommand{\om}{\omega}
\newcommand{\Om}{\Omega}
\newcommand{\de}{\delta}
\newcommand{\De}{\Delta}
\newcommand{\al}{\alpha}

\newcommand{\F}{\mathcal F}

\newcommand{\la}{\lambda}

\newcommand{\Ex}{\mathbb E\,}
\newcommand{\R}{\mathbb R}
\newcommand{\C}{\mathbb C}
\newcommand{\Z}{\mathbb Z}

\newcommand{\Var}{\operatorname{Var}}

\newcommand{\PR}[1]{\mathbb{P}\mskip1.5mu\bigg(\mskip1.5mu#1\mskip1.5mu\bigg)}

\newcommand{\myatop}[2]{{\genfrac{}{}{0pt}{}{#1}{#2}}}

\newcommand{\ABclose}{$(A,B)$\nobreakdash-\hspace{0pt}close}

\newcommand{\Cduplication}{$C$\nobreakdash-\hspace{0pt}duplication}

\newcommand{\splittable}[1]{$#1$\nobreakdash-\hspace{0pt}splittable}
\newcommand{\close}[1]{$#1$\nobreakdash-\hspace{0pt}close}
\newcommand{\integrable}[1]{$#1$\nobreakdash-\hspace{0pt}integrable}

\newcommand{\dimensional}[1]{$#1$\nobreakdash-\hspace{0pt}dimensional}

%% file: intro.tex
Recent results on moderate deviations for random complex zeroes
\cite{NSV} are a challenge for probability theory, since they involve
complex analysis, in contrast to asymptotic normality obtained via random
fields \cite{SoTs1}. Taking up the challenge, I deduce moderate deviations for
random complex zeroes from a new general theorem on moderate deviations for
random fields. As a by-product, the same general theorem gives the asymptotic
normality, avoiding the diagram techniques of \cite{SoTs1}. However, I
consider only smooth test functions, leaving aside interesting effects of
sharp boundary \cite{NSV}.

The random complex zeroes are singularities of a stationary random field on
the complex plane (the logarithm of the absolute value of a normalized
Gaussian analytic function). This random field transcends the existing theory
of moderate deviations (see \cite{DMPU}, \cite{DGW} and references therein) in
several aspects:

(a) it is a random field on the plane, not a random process on the line;

(b) it has some, but not all exponential moments;

(c) it is some function of a Gaussian random field, but the function ($ z
\mapsto \ln |z| $) is singular at $ 0 $, which entangles moderate deviations
of the non-Gaussian field with small deviations of the underlying Gaussian
field;

(d) the underlying Gaussian field is non-stationary.

The first part (Sections \ref{sect1}--\ref{sect4}) contains the general result
in dimension one (random processes on the line). Dimension two (random fields
on the plane) is treated in the second part (Sections
\ref{sect5}--\ref{sect6}), using the first part. The third part (Section
\ref{sect7}) deals with random complex zeros.

Main results of the first part and the third part are formulated below.

\begin{definition}\label{defI1}
A stationary random process $ X = (X_t)_{t\in\R} $ is \emph{splittable,} if
$ \Ex \exp |X_0| < \infty $, $ \Ex X_0 = 0 $, and there exists (on some
probability space) a triple of random processes $ X^0, X^-, X^+ $ such that

(a) the two processes $ X^-, X^+ $ are independent;

(b) the four processes $ X, X^0, X^-, X^+ $ are identically distributed;

(c) $ \Ex \exp \( \int_{-\infty}^0 |X^-_t - X^0_t| \, \D t + \int_0^\infty
|X^+_t - X^0_t| \, \D t \) < \infty $.
\end{definition}

\begin{theorem}\label{theorem1}
For every splittable stationary random process $ X $ there exists $ \si \in
[0,\infty) $ such that for every compactly supported continuous function $ f :
\R \to \R $,
\[
\lim_\myatop{ r\to\infty, \la\to0 }{ \la\log r\to0 }
\frac1{ r\la^2 } \ln \Ex \exp \la \int_{-\infty}^\infty f \Big( \frac t r
\Big) X_t \, \D t = \frac{\si^2}2 \| f \|^2_{L_2(\R)} \, .
\]
\end{theorem}

That is, for every $ \eps $ there exist $ R $ and $ \de $ such that the
given expression is \close{\eps} to the right-hand side for all $ r \ge R $
and all $ \la \ne 0 $ such that $ |\la| \log r \le \de $.

\begin{corollary}\label{corollary3}
Let $ X $, $ \si $ and $ f $ be as in Theorem \ref{theorem1}, and $ \si \ne 0
$. Then
\[
\lim_\myatop{ r\to\infty, c\to\infty }{ (c\log r)^2/r \to0 }
\frac1{c^2} \ln \PR{ \int f \Big( \frac t r \Big) X_t \, \D t \ge c\si \| f
\|_{L_2(\R)} \sqrt r } = -\frac12 \, .
\]
\end{corollary}

Unfortunately, the region of moderate deviations ($ r \to \infty $, $ c \to
\infty $, $ \frac{c^2}r \to 0 $) is not covered. The condition $ \frac{ (c
\log r)^2 }{ r } \to 0 $ leaves a small gap between Corollary \ref{corollary3}
and large deviations ($ \frac{ c^2 } r = \const $).

\begin{corollary}\label{corollary4}
Let $ X $, $ \si $ and $ f $ be as in Theorem \ref{theorem1}. Then the
distribution of $ r^{-1/2} \int f \( \frac t r \) X_t \, \D t $ converges (as
$ r \to \infty $) to the normal distribution $ N(0,\si^2 \| f \|^2) $.
\end{corollary}

Consider now the random entire function $ \psi : \C \to \C $ defined by
\[
\psi(z) = \sum_{k=0}^\infty \frac{ \zeta_k z^k }{ \sqrt{k!} } \, ,
\]
where $ \zeta_0, \zeta_1, \dots $ are independent standard complex Gaussian
random variables.

\begin{theorem}\label{theorem2}
There exists an absolute constant $ \si \in (0,\infty) $ such that for every
compactly supported $ C^2 $-function $ h : \C \to \R $,
\[
\lim_\myatop{ r\to\infty, \la\to0 }{ \la \log^2 r \to 0 }
\frac1{ r^2 \la^2 } \ln \Ex \exp \la r^2 \bigg( \sum_{z:\psi(z)=0} h \Big(
\frac z r \Big) - \frac{ r^2 }\pi { \textstyle \int h \, \D m } \bigg) =
\frac{ \si^2 }2 { \textstyle \int | \De h |^2 \, \D m } \, ;
\]
here $ m $ is the Lebesgue measure on $ \C $, and $ \De h $ is the Laplacian
of $ h $.
\end{theorem}

\begin{corollary}\label{corollary6}
Let $ \si $ and $ h $ be as in Theorem \ref{theorem2}. Then
\[
\hspace*{-5mm}
\lim_\myatop{ r\to\infty, c\to\infty }{ (c\log^2 r)/r \to0 }
\! \frac1{c^2} \ln \PR{ \! \sum_{z:\psi(z)=0} \! h \Big( \frac z r \Big) -
\frac{ r^2 }\pi { \textstyle \int h \, \D m } \ge \frac{ c\si }{ r } \sqrt{
{ \textstyle \int | \De h |^2 \, \D m } } \, } = -\frac12 \, .
\]
\end{corollary}

The same holds for $ (-h) $, of course.

\begin{corollary}\label{corollary7}
Let $ \si $ and $ h $ be as in Theorem \ref{theorem2}. Then the distribution
of
\[
r \, \bigg( \sum_{z:\psi(z)=0} h \Big( \frac z r \Big) - \frac{ r^2 }\pi
{ \textstyle \int h \, \D m } \bigg)
\]
converges (as $ r \to \infty $) to the normal distribution $ N ( 0,\si^2 \int
| \De h |^2 \, \D m ) $.
\end{corollary}

This is the asymptotic normality established in \cite{SoTs1} using moments and
diagrams.

%% file: sect1.tex
Main results of this section, formulated below, are Theorem \ref{theoremA}
(used in Sect.\ \ref{sect2}), Proposition \ref{lemma4} (also used in Sect.\
\ref{sect2}), and Proposition \ref{lemma44} (used in Sect.\ \ref{sect5}, see
\ref{lemma31}).

\begin{definition}\label{definitionA}
Let $ X,Y $ be random variables (possibly on different probability spaces) and
$ C \in [0,\infty) $. We say that $ Y $ is a \emph{\Cduplication} of $ X $, if
there exist random variables $ X_1, X_2, Z $ (on some probability space) such
that

$ X_1 $ and $ X_2 $ are independent,

$ X_1, X_2, X $ are identically distributed,

$ X_1 + X_2 + Z $ and $ Y $ are identically distributed,

$ \ln \Ex \exp \la Z \le C \la^2 $ for all $ \la \in [-1,1] $.
\end{definition}

\begin{theorem}\label{theoremA}
Let random variables $ X_1,X_2,\dots $ (possibly on different probability
spaces) and numbers $ C_1,C_2,\dots \in [0,\infty) $ be such that

(a) $ X_{n+1} $ is a $ C_n $-duplication of $ X_n $ (for all $ n = 1,2,\dots
$);

(b) $ \sup_n \( (2\theta)^{-n} C_n \) < \infty $ for some $ \theta < 1 $;

(c) $ \Ex \exp \eps |X_1| < \infty $ for some $ \eps > 0 $;

(d) $ \Ex X_1 = 0 $.

\noindent Then the following limit exists:
\[
\lim_{n\to\infty, \la n \to 0} \frac1{ 2^n \la^2 } \ln \Ex \exp \la X_n \, .
\]
\end{theorem}

\medskip

Given a function $ f : (0,\infty) \to [0,\infty] $ and a number $ C \in
[0,\infty) $, we define another function $ f_+[C] : (0,\infty) \to
[0,\infty] $ as follows:
\begin{equation}\label{A*}
\begin{aligned}
f_+[C](\la) & = \inf_{p\in[1/(1-\la),\infty)} \; \frac 2 p f(p\la) +
 \frac{p}{p-1} C \la^2 \quad & \text{for } \la \in (0,1) \, , \\
f_+[C](\la) & = \infty \quad & \text{for } \la \in [1,\infty) \, .
\end{aligned}
\end{equation}
Further, we define recursively for $ n=0,1,2,\dots $
\begin{equation}\label{1.4}
f_+[C_0,\dots,C_{n+1}] = ( f_+[C_0,\dots,C_n])_+[C_{n+1}] \, .
\end{equation}

\begin{proposition}\label{lemma4}
For every $ \eps, \theta \in (0,1) $,
\begin{multline*}
\limsup_{n\to\infty,\la n\to0+} \frac1{ 2^{n+1} \la^2 } f_+[C_0,\dots,C_n]
 (\la) \le \\
\le \frac1{1-\eps} \limsup_{\la\to0} \frac1{\la^2} f(\la) + \frac1{ 2\eps
  (1-\sqrt\theta)^2 } \sup_n \frac{ C_n }{ (2\theta)^n } \, .
\end{multline*}
\end{proposition}

\begin{proposition}\label{lemma44}
For every $ \eps, \theta \in (0,1) $, $ n $ and $ \la $,

(a) if $ 0 < \la \le \eps \theta^{n/2} (1-\sqrt\theta) $ then
\[
\frac1{ 2^{n+1} } f_+[C_0,\dots,C_n] (\la) \le (1-\eps) f \Big(
\frac{\la}{1-\eps} \Big) + \frac{\la^2}{2\eps} \frac1{ (1-\sqrt\theta)^2 }
\max_{k=0,\dots,n} \frac{ C_k }{ (2\theta)^k } \, ;
\]

(b) let $ m \in \{0,1,\dots,n-1\} $ be such that $ (n-m) \la < 1 $ and $
\mu \le \eps \theta^{m/2} (1-\sqrt\theta) $, where $ \mu = \frac{ \la }{ 1 -
(n-m) \la } $; then
\begin{multline*}
\frac1{ 2^{n+1} } f_+[C_0,\dots,C_n] (\la) \le \\
\le (1-\eps) f \Big(
\frac{\mu}{1-\eps} \Big) + \bigg( \frac{\mu^2}{2\eps} \frac1{
(1-\sqrt\theta)^2 } + \frac\la2 \frac{ \theta^{m+1} }{ 1-\theta } \bigg)
\max_{k=0,\dots,n} \frac{ C_k }{ (2\theta)^k } \, .
\end{multline*}
\end{proposition}

The proof of Th.~\ref{theoremA} uses Prop.~\ref{lemma4} (whose proof uses
Prop.~\ref{lemma44}). The relevance of \eqref{A*}-\eqref{1.4} and
\ref{lemma4} to \ref{theoremA} stems from Lemmas \ref{lemmaA0}, \ref{lemma1}.
The proofs of Propositions \ref{lemma44} and \ref{lemma4} are given after
Lemmas \ref{lemma2}, \ref{lemma3}. For the proof of Th.~\ref{theoremA} see the
end of this section.

\begin{lemma}\label{lemmaA0}
For all random variables $ X,Y $ and all $ p \in (1,\infty) $,
\begin{multline*}
p \ln \Ex \exp \frac1p X - (p-1) \ln \Ex \exp \Big( - \frac1{p-1} Y \Big) \le
 \\
\le \ln \Ex \exp (X+Y) \le \frac1p \ln \Ex \exp pX + \frac{p-1}p \ln \Ex \exp
\frac{p}{p-1} Y \, .
\end{multline*}
(In the lower bound we interpret $ \infty-\infty $ as $ -\infty $.)
\end{lemma}

\begin{proof}
By the H\"older inequality,
\[
\Ex \exp(X+Y) = \Ex ( \exp X \cdot \exp Y ) \le ( \Ex \exp pX )^{1/p} \Big(
\Ex \exp \frac{p}{p-1} Y \Big)^{(p-1)/p} \, ;
\]
the upper bound follows. We apply the upper bound to $ \frac1p (X+Y) $ and $
\( -\frac1p Y \) $ instead of $ X,Y $:
\[
\ln \Ex \exp \frac1p X \le \frac1p \ln \Ex \exp (X+Y) + \frac{p-1}p \ln \Ex
\exp \Big( -\frac1{p-1} Y \Big) \, ;
\]
the lower bound follows.
\end{proof}

\begin{lemma}\label{lemma1}
(a) Let random variables $ X,Y $ and a function $ f : (0,\infty) \to
[0,\infty] $ be such that
\begin{gather*}
Y \text{ is a \Cduplication\ of } X \, , \\
\ln \Ex \exp \la X \le f(\la) \quad \text{for all } \la \in (0,\infty) \, .
\end{gather*}
Then
\[
\ln \Ex \exp \la Y \le f_+[C](\la) \quad \text{for all } \la \in (0,\infty) \,
.
\]

(b) Let random variables $ X,Y,Z $ and a function $ f : (0,\infty) \to
[0,\infty] $ be such that
\begin{gather*}
X,Y \text{ are independent,} \\
\ln \Ex \exp \la X \le f(\la) \quad \text{and} \quad \ln \Ex \exp \la Y \le
 f(\la) \quad \text{for all } \la \in (0,\infty) \, , \\
\ln \Ex \exp \la Z \le C \la^2 \quad \text{for all } \la \in [0,1] \, .
\end{gather*}
Then
\[
\ln \Ex \exp \la(X+Y+Z) \le f_+[C](\la) \quad \text{for all } \la \in
(0,\infty) \, .
\]
\end{lemma}

\begin{proof}
Item (a) is a special case of Item (b). Item (b) follows from Lemma
\ref{lemmaA0} (the upper bound) applied to $ \la(X+Y) $ and $ \la Z $.
\end{proof}

\begin{lemma}\label{lemma2}
\begin{multline*}
\frac1{2^{n+1}} f_+[C_0,\dots,C_n](\la) \le \( 1 - (n+1)\la \) f \Big(
 \frac{\la}{ 1 - (n+1)\la } \Big) + \frac\la2 \sum_{k=0}^n 2^{-k} C_k \\
\text{for } 0 < \la < \frac1{n+1} \, .
\end{multline*}
\end{lemma}

\begin{proof}
Induction in $ n $. For $ n=0 $, the needed inequality
\begin{equation}\label{A3}
\frac12 f_+[C](\la) \le ( 1 - \la ) f \Big( \frac{\la}{ 1 - \la } \Big) +
\frac\la2 C \quad \text{for } 0<\la<1
\end{equation}
follows from \eqref{A*} for $ p = \frac1{1-\la} $. For $ n>0 $,
denoting $ f_+[C_0,\dots,C_{n-1}] $ by $ g $, the assumed inequality for $ n-1
$ takes the form
\begin{equation}\label{A4}
\frac1{2^{n}} g(\la) \le \( 1 - n\la \) f \Big( \frac{\la}{ 1 - n\la } \Big) +
\frac\la2 \sum_{k=0}^{n-1} 2^{-k} C_k \quad \text{for } 0<\la<\frac1n \, , 
\end{equation}
while the needed inequality for $ n $ becomes
\begin{multline*}
\frac1{2^{n+1}} g_+[C_n](\la) \le \( 1 - (n+1)\la \) f \Big(
 \frac{\la}{ 1 - (n+1)\la } \Big) + \frac\la2 \sum_{k=0}^n 2^{-k} C_k \\
\text{for } 0<\la<\frac1{n+1} \, . 
\end{multline*}

Let $ 0<\la<\frac1{n+1} $. By \eqref{A3}, $ \frac12 g_+[C_n](\la) \le ( 1 - \la
) g ( \mu ) + \frac\la2 C_n $, where $ \mu = \frac{\la}{ 1 - \la } $. We note
that $ 1-n\mu = \frac{ 1 - (n+1)\la }{ 1-\la } > 0 $, $ \frac{\mu}{1-n\mu} =
\frac{\la}{ 1 - (n+1)\la } $ and get from \eqref{A4}
\[
\frac1{2^{n}} g(\mu) \le \frac{ 1 - (n+1)\la }{ 1-\la } f \Big( \frac{\la}{ 1
  - (n+1)\la } \Big) + \frac\mu2 \sum_{k=0}^{n-1} 2^{-k} C_k \, .
\]
Thus,
\begin{multline*}
\frac1{2^{n+1}} g_+[C_n](\la) \le \frac{ 1 - \la }{ 2^n } g (\mu)
 + \frac\la{2^{n+1}} C_n \le \\
\le \( 1 - (n+1)\la \) f \Big( \frac{\la}{ 1 - (n+1)\la } \Big) + \frac\la2
\sum_{k=0}^{n-1} 2^{-k} C_k + \frac\la2 2^{-n} C_n \, .
\end{multline*}
\end{proof}

\begin{lemma}\label{lemma3}
For every $ \eps \in (0,1) $,
\begin{multline*}
\frac1{2^{n+1}} f_+[C_0,\dots,C_n](\la) \le (1-\eps) f \Big(
 \frac{\la}{1-\eps} \Big) + \frac{\la^2}{2\eps} \bigg( \sum_{k=0}^n 2^{-k/2}
 \sqrt{C_k} \bigg)^2 \\
\text{for } 0 < \la \le \eps \frac{ \min_{k=0,\dots,n} 2^{-k/2} \sqrt{C_k} }{
 \sum_{k=0,\dots,n} 2^{-k/2} \sqrt{C_k} } \, .
\end{multline*}
\end{lemma}

\begin{proof}
Induction in $ n $. For $ n=0 $, the needed inequality
\[
\frac12 f_+[C](\la) \le ( 1 - \eps ) f \Big( \frac{\la}{ 1 - \eps } \Big) +
\frac{\la^2}{2\eps} C \quad \text{for } 0<\la\le\eps
\]
follows from \eqref{A*} for $ p = \frac1{1-\eps} $. For $ n>0 $,
we write the assumed inequality for $ n-1 $ in the form
\begin{multline}\label{A7}
\frac1{2^{n}} g(\la_\old) \le ( 1 - \eps_\old ) f \Big( \frac{\la_\old}{ 1 -
 \eps_\old } \Big) + \frac{\la^2_\old}{2\eps_\old} \si^2_\old \\
\text{whenever } 0 < \la_\old \le \eps_\old \frac{ m_\old }{ \si_\old } \text{
 and } 0 < \eps_\old < 1 \, ;
\end{multline}
here $ g = f_+[C_0,\dots,C_{n-1}] $, $ \si_\old = \sum_{k=0}^{n-1} 2^{-k/2}
\sqrt{C_k} $ and $ m_\old = \linebreak[4] \min_{k=0,\dots,n-1} 2^{-k/2}
\sqrt{C_k} $. We have to prove that
\begin{multline}\label{A8}
\frac1{2^{n+1}} g_+[C_n](\la) \le (1-\eps) f \Big(
 \frac{\la}{1-\eps} \Big) + \frac{\la^2}{2\eps} \si^2 \\
\text{whenever } 0 < \la \le \eps \frac{ m }{ \si } \text{ and } 0 < \eps < 1
\, ;
\end{multline}
here $ \si = \sum_{k=0}^{n} 2^{-k/2} \sqrt{C_k} = \si_\old + 2^{-n/2}
\sqrt{C_n} $ and $ m = \min_{k=0,\dots,n} 2^{-k/2} \sqrt{C_k} = \min ( m_\old,
2^{-n/2} \sqrt{C_n} ) $. 

By \eqref{A*} for $ p = \frac1{1-\eps_\new} $,
\begin{multline*}
g_+[C_n](\la) \le 2 (1-\eps_\new) g \Big( \frac{\la}{1-\eps_\new} \Big) +
 \frac{ C_n \la^2 }{ \eps_\new } \\
\text{whenever } 0 < \la \le \eps_\new \text{ and } 0 < \eps_\new < 1 \, .
\end{multline*}
We combine it with \eqref{A7} for $ \la_\old = \frac{ \la }{ 1 - \eps_\new }
$:
\begin{multline*}
\frac1{2^{n+1}} g_+[C_n](\la) \le \\
(1-\eps_\new) \bigg( (1-\eps_\old) f \Big( \frac{ \la }{ (1-\eps_\new)
 (1-\eps_\old) } \Big) + \frac{ \la^2 }{ 2\eps_\old (1-\eps_\new)^2 }
 \si^2_\old \bigg) + \frac{ C_n \la^2 }{ 2^{n+1} \eps_\new } \, .
\end{multline*}
In order to get \eqref{A8} it remains to find $ \eps_\old, \eps_\new \in
(0,1) $ such that $ (1-\eps_\old)(1-\eps_\new) = 1-\eps $,
\begin{equation}\label{A9}
\frac{ \si^2_\old }{ 2\eps_\old (1-\eps_\new) } + \frac{ C_n }{ 2^{n+1}
  \eps_\new } \le \frac{ \si^2 }{ 2\eps }
\end{equation}
and $ \la \le \eps_\new $, $ \frac{ \la }{ 1-\eps_\new } \le \eps_\old \frac{
m_\old }{ \si_\old } $.

We define
\[
x = \frac{ \si_\old }{ 2^{-n/2} \sqrt{C_n} } \, , \quad \eps_\old = \frac{
\eps x }{ x+1-\eps } \, , \quad \eps_\new = \frac{ \eps }{ x+1 } \, .
\]
Clearly, $ \eps_\old < 1 $ (since $ \eps x < x < x+1-\eps $) and $ \eps_\new <
1 $ (since $ \eps < 1 < x+1 $). Also,
\[
(1-\eps_\old)(1-\eps_\new) = \frac{ x+1-\eps-\eps x }{ x+1-\eps } \cdot \frac{
x+1-\eps }{ x+1 } = 1 - \eps \, .
\]
Taking into account that $ 2^{-n/2} \sqrt{C_n} = \frac{ \si_\old }{ x } $ and
$ \frac{ \si_\old }{ \si } = \frac{ x }{ x+1 } = \frac{ (1-\eps_\new)
\eps_\old }{ \eps } $ we get
\begin{multline*}
\la \le \eps \frac{ m }{ \si } = \frac{ \eps }{ \si } \min \Big( m_\old, \frac{
 \si }{ x+1 } \Big) = \\
= \min \Big( \frac{ \eps m_\old }{ \si }, \frac{ \eps }{ x+1 } \Big) = \min
\Big( m_\old \cdot \frac{ (1-\eps_\new) \eps_\old }{ \si_\old }, \eps_\new
\Big) \, .
\end{multline*}
Finally, we check \eqref{A9}:
\begin{multline*}
\frac{ \si^2_\old }{ 2\eps_\old (1-\eps_\new) } + \frac{ C_n }{ 2^{n+1}
 \eps_\new } = \frac{ \si^2_\old }{ 2 } \cdot \frac{ x+1 }{ \eps x } + \frac1{
 2 \eps_\new } \Big( \frac{ \si_\old }{ x } \Big)^2 =\\
= \frac{ \si^2_\old ( x+1 ) }{ 2 \eps x } + \frac{ (x+1) \si^2_\old }{ 2\eps
 x^2 } = \frac{ \si^2_\old ( x+1 ) }{ 2 \eps x^2 } (x+1) = \\
= \frac1{ 2\eps } \Big( \frac{ \si_\old ( x+1 ) }{ x } \Big)^2 = \frac{ \si^2
}{ 2\eps } \, .
\end{multline*}
\end{proof}

The reader may wonder, how did I found these formulas for $ \eps_\old,
\eps_\new $ in the proof of Lemma \ref{lemma3}. In fact, I have minimized the
left-hand side of \eqref{A9} in $ \eps_\old, \eps_\new $ restricted by $
(1-\eps_\old)(1-\eps_\new) = 1-\eps $. The expression $ \( \sum 2^{-k/2}
\sqrt{C_k} \)^2 $ have appeared afterwards.

\begin{proof}[Proof of Prop.~\textup{\ref{lemma44}.}]
Denote $ \max_{k=0,\dots,n} \frac{ C_k }{ (2\theta)^k } $ by $ M $. Without
loss of generality we assume that $ C_k = M (2\theta)^k $ for all $ k $. Then
\[
\sum_{k=0}^n 2^{-k/2} \sqrt{C_k} \le \frac{ \sqrt M }{ 1 - \sqrt\theta } \quad
\text{and} \quad \min_{k=0,\dots,n} 2^{-k/2} \sqrt{C_k} = \sqrt M
\theta^{n/2} \, .
\]
Item (a) follows immediately from Lemma \ref{lemma3}.

Proving Item (b), we use Lemma \ref{lemma3} for passing from $ f $ to $ g =
\linebreak
f_+[C_0,\dots,C_m] $, and Lemma \ref{lemma2} for passing from $ g $ to
$ g_+[C_{m+1},\dots,C_n] = f_+[C_0,\dots,C_n] $. Namely,
\[
\frac1{2^{m+1}} f_+[C_0,\dots,C_m](\mu) \le (1-\eps) f \Big(
\frac{\mu}{1-\eps} \Big) + \frac{\mu^2}{2\eps} \frac{ M }{ (1-\sqrt\theta)^2 }
\]
and
\begin{multline*}
\frac1{ 2^{n-m} } g_+[C_{m+1},\dots,C_n] (\la) \le \\
\le \( 1 - (n-m) \la \) g \Big( \frac{ \la }{ 1 - (n-m) \la } \Big) +
 \frac{\la}{2} \sum_{k=m+1}^n 2^{-(k-m-1)} C_k \, .
\end{multline*}
\end{proof}

\begin{proof}[Proof of Prop.~\textup{\ref{lemma4}.}]
Let $ \sup_n \frac{ C_n }{ (2\theta)^n } = M < \infty $ (otherwise there is
nothing to prove). We consider two overlapping cases separately.

The first case: $ \la \ll \theta^{n/2} $.

Prop.\ \ref{lemma44}(a) gives eventually (when $ \theta^{-n/2} \la $ is small
enough),
\[
\frac1{2^{n+1}} f_+[C_0,\dots,C_n](\la) \le (1-\eps) f \Big(
\frac{\la}{1-\eps} \Big) + \frac{\la^2}{2\eps} \frac{ M }{ (1-\sqrt\theta)^2 }
\, .
\]
It remains to divide by $ \la^2 $ and note that
\[
(1-\eps) \limsup_{\la\to0+} \frac1{ \la^2 } f \Big( \frac{\la}{1-\eps} \Big) =
\frac1{ 1-\eps } \limsup_{\la\to0+} \frac1{ \la^2 } f(\la) \, .
\]

The second case: $ \la \gg \theta^{n} $.

We choose $ m \in \{ 0,1,\dots,n-1 \} $ such that
\[
\la^2 \ll \theta^{m} \ll \la \, ,
\]
apply Prop.\ \ref{lemma44}(b), observe that $ \mu \sim \la $, $ \frac\la2
\frac{ \theta^{m+1} }{ 1-\theta } = o(\la^2) $ and get
\[
\frac1{2^{n+1}} f_+[C_0,\dots,C_n](\la) \le (1-\eps) f \Big(
\frac{\mu}{1-\eps} \Big) + \bigg( \frac{\mu^2}{2\eps} \frac1{
  (1-\sqrt\theta)^2 } + o \( \la^2 \) \bigg) M \, .
\]
It remains to divide by $ \la^2 $ and note that
\[
(1-\eps) \limsup \frac1{ \la^2 } f \Big( \frac{\mu}{1-\eps} \Big) =
\frac1{ 1-\eps } \limsup_{\la\to0+} \frac1{ \la^2 } f(\la) \, .
\]
\end{proof}

We are not yet in position to prove Theorem \ref{theoremA}, since the chains
of inequalities considered above are based on the upper bound of Lemma
\ref{lemmaA0}. We need similar results on lower bounds.

Given $ f : (0,\infty) \to [0,\infty] $ and $ C \in [0,\infty) $, we define $
f_-[C] : (0,\infty) \to [0,\infty] $ by
\begin{equation}\label{A**}
f_-[C](\la) = \sup_{p\in[\la+1,\infty)} \; 2p f \Big( \frac\la p \Big) -
\frac1{p-1} C \la^2 \quad \text{for } \la \in (0,\infty) \, .
\end{equation}
(The supremum is nonnegative, since the limit as $ p \to \infty $ is
nonnegative.)
Further, we define recursively for $ n=0,1,2,\dots $
\begin{equation}
f_-[C_0,\dots,C_{n+1}] = ( f_-[C_0,\dots,C_n])_-[C_{n+1}] \, .
\end{equation}

Lemmas \ref{lemma11}, \ref{lemma12}, \ref{lemma13} and Propositions
\ref{lemma444}, \ref{lemma14} are lower-bound counterparts of \ref{lemma1},
\ref{lemma2}, \ref{lemma3}, \ref{lemma44} and \ref{lemma4}. Proofs are quite
similar, but many small changes in formulas look unpredictable, especially in
the proof of \ref{lemma13}.

\begin{lemma}\label{lemma11}
(a) Let random variables $ X,Y $ and a function $ f : (0,\infty) \to
[0,\infty] $ be such that
\begin{gather*}
Y \text{ is a \Cduplication\ of } X \, , \\
\ln \Ex \exp \la X \ge f(\la) \quad \text{for all } \la \in (0,\infty) \, .
\end{gather*}
Then
\[
\ln \Ex \exp \la Y \ge f_-[C](\la) \quad \text{for all } \la \in (0,\infty) \,
.
\]

(b) Let random variables $ X,Y,Z $ and a function $ f : (0,\infty) \to
[0,\infty] $ be such that
\begin{gather*}
X,Y \text{ are independent,} \\
\ln \Ex \exp \la X \ge f(\la) \quad \text{and} \quad \ln \Ex \exp \la Y \ge
 f(\la) \quad \text{for all } \la \in (0,\infty) \, , \\
\ln \Ex \exp \la Z \le C \la^2 \quad \text{for all } \la \in [0,1] \, .
\end{gather*}
Then
\[
\ln \Ex \exp \la(X+Y-Z) \ge f_-[C](\la) \quad \text{for all } \la \in
(0,\infty) \, .
\]
\end{lemma}

\begin{proof}
Item (a) is a special case of Item (b). Item (b) follows from Lemma
\ref{lemmaA0} (the lower bound) applied to $ \la(X+Y) $ and $ (-\la Z) $.
\end{proof}

\begin{lemma}\label{lemma12}
\begin{multline*}
\frac1{2^{n+1}} f_-[C_0,\dots,C_n](\la) \ge \( 1 + (n+1)\la \) f \Big(
 \frac{\la}{ 1 + (n+1)\la } \Big) - \frac\la2 \sum_{k=0}^n 2^{-k} C_k \\
\text{for } 0 < \la < \infty \, .
\end{multline*}
\end{lemma}

\begin{proof}
Induction in $ n $. For $ n=0 $, the needed inequality
\begin{equation}\label{A16}
\frac12 f_-[C](\la) \ge ( 1 + \la ) f \Big( \frac{\la}{ 1 + \la } \Big) -
\frac\la2 C \quad \text{for } 0<\la<\infty
\end{equation}
follows from \eqref{A**} for $ p = 1+\la $. For $ n>0 $,
denoting $ f_-[C_0,\dots,C_{n-1}] $ by $ g $, the assumed inequality for $ n-1
$ takes the form
\begin{equation}\label{A17}
\frac1{2^{n}} g(\la) \ge \( 1 + n\la \) f \Big( \frac{\la}{ 1 + n\la } \Big) -
\frac\la2 \sum_{k=0}^{n-1} 2^{-k} C_k \quad \text{for } 0<\la<\infty \, , 
\end{equation}
while the needed inequality for $ n $ becomes
\[
\frac1{2^{n+1}} g_-[C_n](\la) \ge \( 1 + (n+1)\la \) f \Big(
\frac{\la}{ 1 + (n+1)\la } \Big) - \frac\la2 \sum_{k=0}^n 2^{-k} C_k
\;\; \text{for } 0<\la<\infty \, . 
\]

By \eqref{A16}, $ \frac12 g_-[C_n](\la) \ge ( 1 + \la ) g ( \mu ) - \frac\la2 C_n
$, where $ \mu = \frac{\la}{ 1 + \la } $. We note that $ 1+n\mu = \frac{ 1 +
(n+1)\la }{ 1+\la } $, $ \frac{\mu}{1+n\mu} = \frac{\la}{ 1 + (n+1)\la } $ and
get from \eqref{A17}
\[
\frac1{2^{n}} g(\mu) \ge \frac{ 1 + (n+1)\la }{ 1+\la } f \Big( \frac{\la}{ 1
+ (n+1)\la } \Big) - \frac\mu2 \sum_{k=0}^{n-1} 2^{-k} C_k \, .
\]
Thus,
\begin{multline*}
\frac1{2^{n+1}} g_-[C_n](\la) \ge \frac{ 1 + \la }{ 2^n } g (\mu)
 - \frac\la{2^{n+1}} C_n \ge \\
\ge \( 1 + (n+1)\la \) f \Big( \frac{\la}{ 1 + (n+1)\la } \Big) - \frac\la2
\sum_{k=0}^{n-1} 2^{-k} C_k - \frac\la2 2^{-n} C_n \, .
\end{multline*}
\end{proof}

\begin{lemma}\label{lemma13}
For every $ \eps \in (0,1) $,
\begin{multline*}
\frac1{2^{n+1}} f_-[C_0,\dots,C_n](\la) \ge \frac1{1-\eps} f ( (1-\eps) \la )
 - \frac{1-\eps}{2\eps} \la^2 \bigg( \sum_{k=0}^n 2^{-k/2} \sqrt{C_k} \bigg)^2
 \\
\text{for } 0 < \la \le \frac\eps{1-\eps} \frac{ \min_{k=0,\dots,n} 2^{-k/2}
 \sqrt{C_k} }{ \sum_{k=0,\dots,n} 2^{-k/2} \sqrt{C_k} } \, .
\end{multline*}
\end{lemma}

\begin{proof}
Induction in $ n $. For $ n=0 $, the needed inequality
\[
\frac12 f_-[C](\la) \ge \frac1{1-\eps} f ( (1-\eps) \la ) -
\frac{1-\eps}{2\eps} \la^2 C \quad \text{for } 0<\la\le\frac\eps{1-\eps}
\]
follows from \eqref{A**} for $ p = \frac1{1-\eps} $. For $ n>0 $,
we write the assumed inequality for $ n-1 $ in the form
\begin{multline}\label{A20}
\frac1{2^{n}} g(\la_\old) \ge \frac1{ 1-\eps_\old } f ( (1 - \eps_\old)
 \la_\old ) - \frac{1-\eps_\old}{2\eps_\old} \la^2_\old \si^2_\old \\
\text{whenever } 0 < \la_\old \le \frac{\eps_\old}{1-\eps_\old} \frac{ m_\old
 }{ \si_\old } \text{ and } 0 < \eps_\old < 1 \, ;
\end{multline}
here $ g = f_-[C_0,\dots,C_{n-1}] $, $ \si_\old = \sum_{k=0}^{n-1} 2^{-k/2}
\sqrt{C_k} $ and $ m_\old = \linebreak[4]
\min_{k=0,\dots,n-1} 2^{-k/2} \sqrt{C_k} $. We have to prove that
\begin{multline}\label{A21}
\frac1{2^{n+1}} g_-[C_n](\la) \ge \frac1{1-\eps} f ( (1-\eps) \la ) -
 \frac{1-\eps}{2\eps} \la^2 \si^2 \\
\text{whenever } 0 < \la \le \frac\eps{1-\eps} \frac m \si \text{ and } 0 <
\eps < 1 \, ;
\end{multline}
here $ \si = \sum_{k=0}^{n} 2^{-k/2} \sqrt{C_k} = \si_\old + 2^{-n/2}
\sqrt{C_n} $ and $ m = \min_{k=0,\dots,n} 2^{-k/2} \sqrt{C_k} = \min ( m_\old,
2^{-n/2} \sqrt{C_n} ) $. 

By \eqref{A**} for $ p = \frac1{1-\eps_\new} $,
\begin{multline*}
g_-[C_n](\la) \ge \frac2{1-\eps_\new} g ( (1-\eps_\new) \la ) -
 \frac{ 1-\eps_\new }{ \eps_\new } C_n \la^2 \\
\text{whenever } 0 < \la \le \frac{\eps_\new}{1-\eps_\new} \text{ and } 0 <
\eps_\new < 1 \, .
\end{multline*}
We combine it with \eqref{A20} for $ \la_\old = (1-\eps_\new) \la $:
\begin{multline*}
\frac1{2^{n+1}} g_-[C_n](\la) \ge \frac1{1-\eps_\new} \bigg(
 \frac1{1-\eps_\old} f \( (1-\eps_\new)(1-\eps_\old) \la \) - \\
- \frac{ 1-\eps_\old }{ 2\eps_\old } (1-\eps_\new)^2 \la^2 \si^2_\old \bigg) -
\frac{ 1-\eps_\new }{ 2^{n+1} \eps_\new } C_n \la^2 \, .
\end{multline*}
In order to get \eqref{A21} it remains to find $ \eps_\old, \eps_\new \in
(0,1) $ such that $ (1-\eps_\old)(1-\eps_\new) = 1-\eps $,
\begin{equation}\label{A22}
\frac{ 1-\eps }{ 2\eps_\old } \si^2_\old + \frac{ 1-\eps_\new }{ 2^{n+1}
  \eps_\new } C_n \le \frac{ 1-\eps }{ 2\eps } \si^2
\end{equation}
and $ \la \le \frac{\eps_\new}{1-\eps_\new} $, $ (1-\eps_\new) \la \le
\frac{\eps_\old}{1-\eps_\old} \frac{ m_\old }{ \si_\old } $.

We define
\[
x = \frac{ \si_\old }{ 2^{-n/2} \sqrt{C_n} } \, , \quad \eps_\old = \frac{
\eps x }{ x+1 } \, , \quad \eps_\new = \frac{ \eps }{ x+1-\eps x } \, .
\]
Clearly, $ \eps_\old < 1 $ (since $ \eps x < x < x+1 $) and $ 0 < \eps_\new <
1 $ (since $ \eps + \eps x < x+1 $). Also,
\[
(1-\eps_\old)(1-\eps_\new) = \frac{ x+1-\eps x }{ x+1 } \cdot \frac{
x+1-\eps x - \eps }{ x+1-\eps x } = 1 - \eps \, .
\]
Taking into account that $ 2^{-n/2} \sqrt{C_n} = \frac{ \si_\old }{ x } $,
$ \frac{ \si_\old }{ \si } = \frac{ x }{ x+1 } $ and $ \frac{1-\eps}\eps
\frac{ \eps_\new }{ 1-\eps_\new } = \frac1{x+1} $ we get
\begin{multline*}
\la \le \frac\eps{1-\eps} \frac{ m }{ \si } = \frac{ \eps }{ 1-\eps } \min
 \Big( m_\old, \frac{ \si_\old }{ x } \Big) \cdot \frac{x}{x+1}
 \frac1{\si_\old} = \\
= \min \Big( \frac{ \eps x m_\old }{ (1-\eps) (x+1) \si_\old }, \frac{ \eps }{
 (1-\eps)(x+1) } \Big) = \\
= \min \Big( \frac{ \eps_\old }{ (1-\eps_\old) (1-\eps_\new) }
 \frac{m_\old}{\si_\old}, \frac{ \eps_\new }{ 1-\eps_\new } \Big) \, .
\end{multline*}
Finally, we check \eqref{A22}:
\begin{multline*}
\frac{ 1-\eps }{ 2\eps_\old } \si^2_\old + \frac{ 1-\eps_\new }{ 2\eps_\new
 } \cdot \frac{ C_n }{ 2^n } = \frac{ (1-\eps)(x+1) }{ 2\eps x } \si^2_\old 
 + \frac{1-\eps}{2\eps} (x+1) \Big( \frac{ \si_\old }{ x } \Big)^2 = \\
= \frac{ 1-\eps }{ 2\eps } \frac{ x+1 }{ x } \si^2_\old \Big( 1 + \frac1x
\Big) = \frac{ 1-\eps }{ 2\eps } \si^2 \, .
\end{multline*}
\end{proof}

\begin{proposition}\label{lemma444}
For every $ \eps, \theta \in (0,1) $, $ n $ and $ \la $,

(a) if $ 0 < \la \le \eps \theta^{n/2} (1-\sqrt\theta) $ then
\[
\frac1{ 2^{n+1} } f_-[C_0,\dots,C_n] (\la) \ge \frac1{1-\eps} f \(
(1-\eps) \la \) - \frac{\la^2}{2\eps} \frac1{ (1-\sqrt\theta)^2 }
\max_{k=0,\dots,n} \frac{ C_k }{ (2\theta)^k } \, ;
\]

(b) let $ m \in \{0,1,\dots,n-1\} $ be such that $ \mu \le \eps \theta^{m/2}
(1-\sqrt\theta) $, where $ \mu = \frac{ \la }{ 1 + (n-m) \la } $; then
\begin{multline*}
\frac1{ 2^{n+1} } f_-[C_0,\dots,C_n] (\la) \ge \\
\ge \frac1{1-\eps} f \( (1-\eps) \mu \) - \bigg( \frac{\mu^2}{2\eps} \frac1{
(1-\sqrt\theta)^2 } + \frac\la2 \frac{ \theta^{m+1} }{ 1-\theta } \bigg)
\max_{k=0,\dots,n} \frac{ C_k }{ (2\theta)^k } \, .
\end{multline*}
\end{proposition}

\begin{proof}
Denote $ \max_{k=0,\dots,n} \frac{ C_k }{ (2\theta)^k } $ by $ M $. Without
loss of generality we assume that $ C_k = M (2\theta)^k $ for all $ k $. Then
\[
\sum_{k=0}^n 2^{-k/2} \sqrt{C_k} \le \frac{ \sqrt M }{ 1 - \sqrt\theta } \quad
\text{and} \quad \min_{k=0,\dots,n} 2^{-k/2} \sqrt{C_k} = \sqrt M
\theta^{n/2} \, .
\]
Item (a) follows immediately from Lemma \ref{lemma13}.

Proving Item (b), we use Lemma \ref{lemma13} for passing from $ f $ to $ g =
\linebreak
f_-[C_0,\dots,C_m] $, and Lemma \ref{lemma12} for passing from $ g $ to
$ g_-[C_{m+1},\dots,C_n] = f_-[C_0,\dots,C_n] $. Namely,
\[
\frac1{2^{m+1}} f_-[C_0,\dots,C_m](\mu) \ge \frac1{1-\eps} f \( (1-\eps) \mu
\) - \frac{\mu^2}{2\eps} \frac{ M }{ (1-\sqrt\theta)^2 }
\]
and
\begin{multline*}
\frac1{ 2^{n-m} } g_-[C_{m+1},\dots,C_n] (\la) \ge \\
\ge \( 1 + (n-m) \la \) g \Big( \frac{ \la }{ 1 + (n-m) \la } \Big) -
\frac{\la}{2} \sum_{k=m+1}^n 2^{-(k-m-1)} C_k \, .
\end{multline*}
\end{proof}

\begin{proposition}\label{lemma14}
For every $ \eps, \theta \in (0,1) $,
\begin{multline*}
\liminf_{n\to\infty,\la n\to0+} \frac1{ 2^{n+1} \la^2 } f_-[C_0,\dots,C_n]
 (\la) \ge \\
\ge (1-\eps) \liminf_{\la\to0+} \frac1{\la^2} f(\la) - \frac{1-\eps}{ 2\eps
  (1-\sqrt\theta)^2 } \sup_n \frac{ C_n }{ (2\theta)^n } \, .
\end{multline*}
\end{proposition}

\begin{proof}
Let $ \sup_n \frac{ C_n }{ (2\theta)^n } = M < \infty $ (otherwise there is
nothing to prove). We consider two overlapping cases separately.

The first case: $ \la \ll \theta^{n/2} $.

Prop.\ \ref{lemma444}(a) gives eventually (when $ \theta^{-n/2} \la $ is small
enough),
\[
\frac1{2^{n+1}} f_-[C_0,\dots,C_n](\la) \ge \frac1{1-\eps} f \( (1-\eps) \la
\) - \frac{\la^2}{2\eps} \frac{ M }{ (1-\sqrt\theta)^2 } \, .
\]
It remains to divide by $ \la^2 $ and note that
\[
\frac1{1-\eps} \liminf_{\la\to0+} \frac1{ \la^2 } f \( (1-\eps) \la \) =
(1-\eps) \liminf_{\la\to0+} \frac1{ \la^2 } f(\la) \, .
\]

The second case: $ \la \gg \theta^{n} $.

We choose $ m \in \{ 0,1,\dots,n-1 \} $ such that
\[
\la^2 \ll \theta^{m} \ll \la \, ,
\]
apply Prop.\ \ref{lemma444}(b), observe that $ \mu \sim \la $, $ \frac\la2
\frac{ \theta^{m+1} }{ 1-\theta } = o(\la^2) $ and get
\[
\frac1{2^{n+1}} f_-[C_0,\dots,C_n](\la) \ge \frac1{1-\eps} f \(
(1-\eps) \mu \) - \bigg( \frac{\mu^2}{2\eps} \frac1{
  (1-\sqrt\theta)^2 } + o \( \la^2 \) \bigg) M \, .
\]
It remains to divide by $ \la^2 $ and note that
\[
\frac1{1-\eps} \liminf \frac1{ \la^2 } f \( (1-\eps) \mu \) =
(1-\eps) \liminf_{\la\to0+} \frac1{ \la^2 } f(\la) \, .
\]
\end{proof}

Now we combine lower and upper bounds. (See also \ref{lemma311} for a more
general statement.)

\begin{proposition}\label{prop1}
Let numbers $ C_1,C_2,\dots $ satisfy $ \sup_n \( (2\theta)^{-n} C_n \) <
\infty $ for some $ \theta < 1 $, and functions $ f_1,f_2,\dots : (0,\infty)
\to [0,\infty] $ satisfy $
(f_n)_-[C_n] \le f_{n+1} \le (f_n)_+[C_n] $ for all $ n $. Assume existence of
the limit $ \lim_{\la\to0+} \( \la^{-2} f_n(\la) \) \in [0,\infty) $ for every
$ n $. Then the following limit exists:
\[
\lim_{n\to\infty,\la n\to0+} \frac{ f_n(\la) }{ 2^n \la^2 } \in [0,\infty) \,
.
\]
\end{proposition}

\begin{proof}
By Prop.\ \ref{lemma4}, for every $ m $ and $ \eps $,
\begin{multline*}
\limsup_{n\to\infty,\la n\to0+} \frac1{ 2^{n-m} \la^2 }
 (f_m)_+[C_{m+1},\dots,C_n] (\la) \le \\
\le \frac1{ 1-\eps } \lim_{\la\to0+} \frac{ f_m(\la) }{ \la^2 } + \frac1{ 2\eps
 (1-\sqrt\theta)^2 } \sup_{n>m} \frac{ C_n }{ (2\theta)^{n-m-1} } \, .
\end{multline*}
Therefore
\[
\limsup_{n\to\infty,\la n\to0+} \frac{ f_n(\la) }{ 2^n \la^2 }
\le \frac1{ 1-\eps } \lim_{\la\to0+} \frac{ f_m(\la) }{ 2^m \la^2 }
 + \frac{ \theta^{m+1} }{ 4\eps (1-\sqrt\theta)^2 } \sup_n \frac{ C_n }{
 (2\theta)^{n} } < \infty \, .
\]
The last term vanishes as $ m\to\infty $;
\[
\limsup_{n\to\infty,\la n\to0+} \frac{ f_n(\la) }{ 2^n \la^2 }
\le \frac1{ 1-\eps } \liminf_{m\to\infty} \lim_{\la\to0+} \frac{ f_m(\la) }{
 2^m \la^2 } \, .
\]
Now the factor $ 1/(1-\eps) $ may be dropped.

The same argument applies to lower bounds (using Prop.\ \ref{lemma14});
\[
\liminf_{n\to\infty,\la n\to0+} \frac{ f_n(\la) }{ 2^n \la^2 } \ge
\limsup_{m\to\infty} \lim_{\la\to0+} \frac{ f_m(\la) }{ 2^m \la^2 } \, .
\]
It follows that these four numbers are equal.
\end{proof}

\begin{lemma}\label{lemmaA**}
If $ \Ex X = 0 $ then $ \Ex \exp \la X \ge 1 $ for all $ \la \in \R $.
\end{lemma}

\begin{proof}
$ \Ex \E^{\la X} - 1 = \Ex \( \E^{\la X} - 1 - \la X \) \ge 0 $.
\end{proof}

\begin{lemma}\label{lemmaA***}
If $ Y $ is a \Cduplication\ of $ X $ and $ \Ex \exp \eps |X| < \infty $ for
some $ \eps > 0 $, then $ \Ex \exp \de |Y| < \infty $ for some $ \de > 0 $,
and $ \Ex Y = 2 \, \Ex X $.
\end{lemma}

\begin{proof}
We use \ref{lemma1}(a) and note that $ f(\eps) < \infty $ implies $ f_+[C] \(
\frac{\eps}{1+\eps} \) < \infty $ (see \eqref{A3}), therefore $ \Ex \exp
\frac{\eps}{1+\eps} Y < \infty $. The same holds for $ (-X), (-Y) $, thus $
\Ex \exp \frac{\eps}{1+\eps} |Y| < \infty $. For $ Z $ of Def.\
\ref{definitionA} we have $ \Ex Z = \frac{\D}{\D\la} \big|_{\la=0} \ln \Ex
\exp \la Z = 0 $, therefore $ \Ex Y = 2 \, \Ex X $.
\end{proof}

\begin{proof}[Proof of Theorem \textup{\ref{theoremA}}.]
By induction, using Lemma \ref{lemmaA***}, $ \Ex \exp \eps_n |X_n| < \infty $
for some $ \eps_n > 0 $, and $ \Ex X_n = 0 $. By Lemma \ref{lemmaA**}, $ \Ex
\exp \la X_n \ge 1 $ for all $ \la $. By \ref{lemma1}(a) and \ref{lemma11}(a),
the functions
\[
f_n(\la) = \ln \Ex \exp \la X_n
\]
satisfy $ (f_n)_-[C] \le f_{n+1} \le (f_n)_+[C] $. Prop.~\ref{prop1} ensures
existence of the limit
\[
\lim_{n\to\infty, \la n \to 0+} \frac{ f_n(\la) }{ 2^n \la^2 } \, .
\]
The same argument applied to $ (-X_n) $ gives us $ \lim_{n\to\infty, \la n \to
0-} \frac{ f_n(\la) }{ 2^n \la^2 } $. The two limits are equal, since both are
equal to $ \lim_{n\to\infty} \( 2^{-n-1} f''_n(0) \) $.
\end{proof}

%% file: sect2.tex
Theorem \ref{theorem1}, Corollary \ref{corollary3} and Corollary
\ref{corollary4}, formulated in the introduction, are proved in this
section. Theorem \ref{lemmaC133} is used (and generalized) in
Sect.~\ref{sect5} (see \ref{E5}), and Lemma \ref{lemmaC2} is used in Sections
\ref{sect3}, \ref{sect5}. Many arguments of this section (especially, proofs
of \ref{propC}--\ref{lemmaC12}) are reused in Sect.\ \ref{sect5} (especially,
\ref{5.13}--\ref{5.19}).

Splittability is defined in the introduction (see Def.~\ref{defI1}) for
stationary random processes. The reader may restrict himself to random
processes with continuous sample functions, but in full generality, a random
process is treated in this work as a measurable map $ X : \R \times \Om \to \R
$ (given a probability space $ (\Om,\F,P) $). Two such maps $ X_1,X_2 : \R
\times \Om \to \R $ satisfying
\begin{equation}\label{same}
P \( \{ \om : X_1(t,\om) = X_2(t,\om) \} \) = 1 \quad \text{for every } t \in
\R \, ,
\end{equation}
are the same for all purposes of this work. (Accordingly, we may treat a
random process as an equivalence class.) Every $ t \in \R $ leads to a random
variable $ X_t : \Om \to \R $, $ X_t(\om) = X(t,\om) $, treated as an
equivalence class. Sample paths $ t \mapsto X(t,\om) $ are defined for almost
all $ \om $, but only as equivalence classes; their (dis-)continuity is
irrelevant.

Here is some background (mostly for non-probabilists).

Condition \eqref{same} ensures that $ \int_0^1 | X_1(t,\om) | \, \D t =
\int_0^1 | X_2(t,\om) | \, \D t $ for almost all $ \om $. Thus, the random
variable $ \int_0^1 | X_t | \, \D t $ is well-defined as an equivalence
class of functions $ \Om \to [0,\infty] $.

Two processes $ X : \R \times \Om_1 \to \R $, $ Y : \R \times \Om_2 \to \R $
are called \emph{identically distributed,} if for all $ n $ and $
t_1,\dots,t_n \in \R $ the random vectors $ (X_{t_1},\dots,X_{t_n}) $ and $
(Y_{t_1},\dots,Y_{t_n}) $ are identically distributed. It follows that the
random variables $ \int_0^1 | X_t | \, \D t $ and $ \int_0^1 | Y_t | \, \D t $
are identically distributed.

The same holds for many other integrals, of course. Some of them are mentioned
in Def.~\ref{defI1}, Th.~\ref{theorem1} and the corollaries.

A process $ X $ is called \emph{stationary,} if for every $ s \in \R $ the
shifted process $ (t,\om) \mapsto X(s+t,\om) $ is distributed like $ X $.

Two processes $ X,Y : \R \times \Om \to \R $ are called \emph{independent,} if
for all $ n $ and $ t_1,\dots,t_n \in \R $ the random vectors $
(X_{t_1},\dots,X_{t_n}) $ and $ (Y_{t_1},\dots,Y_{t_n}) $ are independent.

If $ X $ is stationary then the distribution of $ \int_s^{s+1} |X_t| \, \D t $
does not depend on $ s $. If $ X $ and $ Y $ are independent then the random
variables $ \int_0^1 | X_t | \, \D t $ and $ \int_0^1 | Y_t | \, \D
t $ are independent.

By the Fubini theorem, $ \Ex \int_0^1 | X_t | \, \D t = \int_0^1 \Ex | X_t |
\, \D t \in [0,\infty] $, and if this value is finite then $ \int_0^1 X_t \,
\D t $ is a well-defined random variable $ \Om \to \R $, and $ \Ex \int_0^1
X_t \, \D t = \int_0^1 \Ex X_t \, \D t $.

We will return to random processes after two quite general lemmas about random variables.

\begin{lemma}\label{lemmaC1}
For every random variable $ X $ and numbers $ C \in [0,\infty) $, $ \la \in
[-1,1] $,
\[
\text{if } \Ex \cosh X \le \cosh C \quad \text{then} \quad \Ex \cosh \la X \le
\cosh \la C \, .
\]
\end{lemma}

\begin{proof}
It is sufficient to find a number $ u>0 $ such that
\[
\cosh \la X - \cosh \la C \le u ( \cosh X - \cosh C ) \quad \text{a.s.,}
\]
which boils down to the inequality
\[
\cosh \la x - \cosh \la y \le \frac{ \la \sinh \la y }{ \sinh y } ( \cosh x -
\cosh y )
\]
for all $ x \in \R $, $ y \in \R \setminus \{0\} $ and $ \la \in [-1,1] $.
In terms of the continuous function $ f : \R \times \R \to \R $ defined by
\begin{align*}
f(x,y) &= \frac{ \cosh \la x - \cosh \la y }{ \cosh x - \cosh y } \quad
 \text{for } x \ne y \, ,  \\
f(x,x) &= \frac{ \la \sinh \la x }{ \sinh x } \quad \text{for } x \ne 0 \, ,
\quad f(0,0) = \la^2 \, ,
\end{align*}
the needed inequality becomes
\[
f(y,y) \le f(x,y) = f(y,x) \le f(x,x) \quad \text{whenever } 0 < x < y \, .
\]
In order to prove the latter we note that $ \frac{\D}{\D x} f(x,x) \le 0 $ for
$ x>0 $ (since $ \la \tanh x \le \tanh \la x $), and
\[
\frac{\pd}{\pd x} f(x,y) = - \frac{ \sinh x }{ \cosh x - \cosh y } \( f(x,y) -
f(x,x) \) \, .
\]
The function $ g(x) = f(x,y)-f(x,x) $ satisfies $ g(y)=0 $ and
\[
g(x)>0 \imply g'(x)>0 \quad \text{for all } x \in (0,y) \, ,
\]
which shows that $ g(x) \le 0 $ for all $ x \in (0,y) $. Similarly,
\[
g(x)<0 \imply g'(x)>0 \quad \text{for all } x \in (y,\infty) \, ,
\]
which shows that $ g(x) \ge 0 $ for all $ x \in (y,\infty) $.
\end{proof}

\begin{lemma}\label{lemmaC2}
For every random variable $ X $ and numbers $ C \in [0,\infty) $, $ \la \in
[-1,1] $,
\[
\text{if } \; \Ex \exp |X| \le \cosh C \; \text{ and } \; \Ex X = 0 \quad
\text{then} \;\quad\; \Ex \exp \la X \le \exp ( A C^2 \la^2 )
\]
where $ A \in (0,\infty) $ is an absolute constant.
\end{lemma}

\begin{proof}
We take
\[
A = \sup_{x\in(0,\infty)} \frac{ \ln ( 2 \cosh x - 1 ) }{ x^2 } \, ;
\]
the supremum is finite, since the fraction tends to $ 1 $ as $ x \to 0+ $ and
to $ 0 $ as $ x \to \infty $. Using Lemma \ref{lemmaC1} and taking into
account that $ \cosh X \le \exp |X| $ we get
\[
2 \, \Ex \cosh \la X \le 2 \cosh \la C \le 1 + \exp ( A(\la C)^2 ) \, .
\]
On the other hand, $ \Ex \exp (-\la X) \ge 1 $ by Lemma \ref{lemmaA**}. Thus,
\[
\Ex \exp \la X = 2 \, \Ex \cosh \la X - \Ex \exp (-\la X) \le 2 \, \Ex \cosh
\la X - 1 \le \exp ( A C^2 \la^2 ) \, .
\]
\end{proof}

\begin{remark}
Lemma \ref{lemmaC2} is a uniform version of equivalence between the following
two conditions on a random variable $ X $:

(a) there exists $ C < \infty $ such that $ \ln \Ex \exp \la X \le C \la^2 $
for all $ \la \in [-1,1] $;

(b) $ \Ex \exp |X| < \infty $ and $ \Ex X = 0 $.

\noindent This equivalence is easy to prove without Lemma \ref{lemmaC1}, as
follows.

(a) \imp\ (b): $ \Ex \exp |X| \le \Ex (\E^X+\E^{-X}) \le 2 \E^C < \infty $;
and $ \Ex X = \frac{\D}{\D\la} \big|_{\la=0} \ln \Ex \exp \la X = 0 $.

(b) \imp\ (a): the function $ \la \mapsto \la^{-2} \ln \Ex \exp \la X $ is
continuous on $ [-1,0) \cup (0,1] $ and has a finite limit at $ 0 $, therefore
it is bounded.
\end{remark}

We return to random processes.

\begin{lemma}\label{easy}
\[
\Ex \exp \frac{\la}{r} \bigg| \int_0^r X_t \, \D t \bigg| \le \frac1r \int_0^r
\Ex \exp \la |X_t| \, \D t \, .
\]
\end{lemma}

\begin{proof}
By convexity of the function $ x \mapsto \E^{\la x} $,
\[
\exp \bigg( \la \frac{1}{r} \int_0^r |X_t| \, \D t \bigg) \le \frac1r \int_0^r
\exp \la |X_t| \, \D t \, ;
\]
we take the expectation.
\end{proof}

Here are the main steps toward Theorem \ref{theorem1}. Using Th.\
\ref{theoremA} we get the convergence for $ f(\cdot) = 1 $ and $ r \in \{
1,2,4,8,\dots \} $ (Prop.\ \ref{propC}), then tor all $ r $ (Prop.\
\ref{propositionC12}). By `concatenation' (Lemma \ref{lemmaC13}) we treat step
functions $ f $. Continuous $ f $ are reached by a limiting procedure
(Lemma \ref{lemmaC14}) based on an upper bound (Th.\ \ref{lemmaC133}) that
holds for all bounded measurable $ f $; its proof uses Prop.\ \ref{lemma4}.

According to Def.\ \ref{defI1}, all splittable processes are stationary in
this section (but not in Sect.\ \ref{sect5}, see Def.\ \ref{defE6}).

\begin{lemma}\label{lemmaC4}
Let $ X $ be splittable. Then there exists $ C < \infty $ such that

(a) for every $ r \in (0,\infty) $, the random variable $ \int_0^{2r} X_t \,
\D t $ is a \Cduplication\ of the random variable $ \int_0^r X_t \, \D t $;

(b) for every $ s,r \in (0,\infty) $ and every measurable $ f : (-s,r) \to
[-1,1] $ there exist (on some probability space) random variables $ Y_-, Y_+ $
and $ Z $ such that

$ Y_-, Y_+ $ are independent;

$ Y_- $ is distributed like $ \int_{-s}^0 f(t) X_t \, \D t $;

$ Y_+ $ is distributed like $ \int_0^r f(t) X_t \, \D t $;

$ Y_- + Y_+ + Z $ is distributed like $ \int_{-s}^r f(t) X_t \, \D t $;

$ \ln \Ex \exp \la Z \le C \la^2 $ for all $ \la \in [-1,1] $.

\medskip

(c) for $ s,r $ and $ f $ as in (b) and $ \la \in [-1,1] $,
\begin{multline}\label{C6a}
\ln \Ex \exp \la \int_{-s}^r f(t) X_t \, \D t \le \frac1p \ln \Ex \exp p \la
 \int_{-s}^0 f(t) X_t \, \D t + \\
+ \frac1p \ln \Ex \exp p \la \int_0^r f(t) X_t \, \D t + \frac{p}{p-1} C \la^2
\quad \text{for } p \in \Big[ \frac1{1-\la},\infty \Big) \, ;
\end{multline}
\vspace{-4mm}
\begin{multline}\label{C6b}
\ln \Ex \exp \la \int_{-s}^r f(t) X_t \, \D t \ge p \ln \Ex \exp \frac{\la}{p}
 \int_{-s}^0 f(t) X_t \, \D t + \\
+ p \ln \Ex \exp \frac{\la}{p} \int_0^r f(t) X_t \, \D t - \frac1{p-1} C \la^2
\quad \text{for } p \in [1+\la,\infty) \, .
\end{multline}
\end{lemma}

\begin{proof}
Item (a) is a special case of (b). Item (c) follows from (b) by
\ref{lemma1}(b) and \ref{lemma11}(b) applied twice: to $ Y_-, Y_+, Z $ and $
(-Y_-), (-Y_+), (-Z) $. In order to prove Item (b) we take $
X^0, X^-, X^+ $ as in Def.~\ref{defI1}. Random variables $ Y_- = \int_{-s}^0
f(t) X_t^- \, \D t $ and $ Y_+ = \int_0^r f(t) X_t^+ \, \D t $ are independent
and distributed like $ \int_{-s}^0 f(t) X_t \, \D t $ and $ \int_0^r f(t) X_t
\, \D t $ respectively. We define $ Z = \int_{-s}^r f(t) X_t^0 \, \D t - Y_- -
Y_+ $, then the random variable $ Y_- + Y_+ + Z $ is distributed like $
\int_{-s}^{r} f(t) X_t \, \D t $. It remains to ensure that $ \ln \Ex \exp \la
Z \le C \la^2 $ for $ \la \in [-1,1] $.

We have $ \Ex \exp |Z| \le \Ex \exp \( \int_{-s}^0 |f(t)| | X_t^0 - X_t^- | \,
\D t + \int_0^r |f(t)| | X_t^0 - X_t^+ | \, \D t \) \linebreak[0] \le B
$, where $ B = \Ex \exp \( \int_{-\infty}^0 | X_t^0 - X_t^- | \, \D t +
\int_0^\infty | X_t^0 - X_t^+ | \, \D t \) $. Also, $ \Ex Z = 0 $ (the use of
Fubini theorem is justified by finiteness of $ B $). By Lemma \ref{lemmaC2}, $
\ln \Ex \exp \la Z \le A M^2 \la^2 $ if $ M $ is defined by $ \cosh M = B $.
\end{proof}

\begin{remark}
The constant $ C $ in \ref{lemmaC4} depends only on the expectation in Def.\
\ref{defI1}(c).
\end{remark}

\begin{proposition}\label{propC}
If a process $ X $ is splittable then for every $ r \in (0,\infty) $ the
following limit exists:
\[
\lim_{n\to\infty,\la n\to 0} \frac1{2^n \la^2} \ln \Ex \exp \la \int_0^{2^n r}
X_t \, \D t \, .
\]
\end{proposition}

\begin{proof}
Random variables $ X_n = \int_0^{2^n r} X_t \, \D t $ satisfy Conditions
(a)--(d) of Theorem \ref{theoremA}. Indeed, Lemma \ref{lemmaC4}(a) verifies
(a) for $ C_n = C $; (b) and (d) are evident; (c) follows from Lemma
\ref{easy}.
\end{proof}

\begin{lemma}\label{lemmaC10}
Let $ X $ be splittable. Then there exists $ \si \in [0,\infty) $ such that
\[
\lim_{n\to\infty,\la n\to0} \frac1{ 2^n \la^2 } \ln \Ex \exp \la \int_0^{2^nr}
X_t \, \D t = \frac12 \si^2 r \quad \text{for all } r \in (0,\infty) \, .
\]
\end{lemma}

\begin{proof}
For each $ r $ the limit exists by Prop.~\ref{propC}. Denote the limit by $
\phi(r) $. Applying Lemma \ref{lemmaC4}(c) to $ 2^n s, 2^n r $ (and $ f(\cdot)
= 1 $), dividing \eqref{C6a}, \eqref{C6b} by $ 2^n \la^2 $ and taking the
limit we get
\[
\frac1p \( \phi(s) + \phi(r) \) \le \phi(s+r) \le p \( \phi(s) + \phi(r) \)
\]
for all $ p > 1 $; thus, $ \phi(s+r) = \phi(s) + \phi(r) $ for all $ s,r > 0
$. Being measurable, such $ \phi $ must be linear. Finally, $ \phi(r) \ge 0 $,
since $ \Ex \exp \la \int_0^r X_t \, \D t \ge 1 $ by Lemma \ref{lemmaA**}.
\end{proof}

\begin{lemma}\label{lemmaC11}
The convergence in Lemma \ref{lemmaC10} is uniform in $ r \in (0,M) $ for
every $ M>0 $.
\end{lemma}

\begin{proof}
Denoting
\[
f_n(\la,r) = \frac1{ 2^n \la^2 } \ln \Ex \exp \la \int_0^{2^n r}
X_t \, \D t - 0.5 \si^2 r \, ,
\]
we have (for every $ t $) $ f_n(\la,r) \to 0 $ as $ n\to\infty $, $ \la n \to
0 $. By \eqref{C6a},
\begin{multline*}
2^n \la^2 \( f_n(\la,s+r) + 0.5 \si^2(s+r) \) \le \frac1p 2^n
 (p\la)^2 \( f_n(p\la,s) + 0.5 \si^2 s \) + \\
+ \frac1p 2^n (p\la)^2 \( f_n(p\la,r) + 0.5 \si^2 r \) +
\frac{p}{p-1} C \la^2
\end{multline*}
provided that $ (1-\la) p \ge 1 $. Thus,
\[
f_n(\la,s+r) \le p f_n(p\la,s) + p f_n(p\la,r) + (p-1) \cdot 0.5 \si^2(s+r) +
2^{-n} \frac{p}{p-1} C \, .
\]
By \eqref{C6b},
\begin{multline*}
2^n \la^2 \( f_n(\la,s+r) + 0.5 \si^2(s+r) \) \ge p 2^n \Big(
 \frac{\la}p \Big)^2 \Big( f_n \Big( \frac{\la}p, s \Big) + 0.5 \si^2 s
 \Big) + \\
+ p 2^n \Big( \frac{\la}p \Big)^2 \Big( f_n \Big( \frac{\la}p, r \Big) +
 0.5 \si^2 r \Big) - \frac1{p-1} C \la^2
\end{multline*}
provided that $ p-\la \ge 1 $. Thus,
\[
f_n(\la,s+r) \ge \frac1p f_n \Big( \frac{\la}p, r \Big) + \frac1p f_n \Big(
\frac{\la}p, r \Big) - \frac{p-1}p \cdot 0.5 \si^2(s+r) - 2^{-n} \frac1{p-1} C
\, .
\]
Given $ \eps > 0 $, we take $ N $ and $ \de > 0 $ such that $ |f_n(\la,r)| \le
\eps $ whenever $ n \ge N $, $ n |\la| \le \de $ and $ r \in A $ for some
measurable set $ A \subset (0,2M) $ of Lebesgue measure $ > \frac32 M $.
Increasing $ N $ if needed, we guarantee that $ \frac{\de}{N} \le \frac12 
$, $ \frac{C}{2^N \eps} \le \frac12 $,
\[
\Big( 1 - \frac{\de}{N} \Big) \Big( 1 + \frac{\eps}{\si^2 M} \Big) \ge 1 \quad
\text{and} \quad \Big( 1 - \frac{C}{2^N \eps} \Big) \Big( 1 +
\frac{\eps}{\si^2 M} \Big) \ge 1 \, .
\]
We take $ p $ such that
\[
\Big( 1 - \frac{\de}{N} \Big) p \ge 1 \, , \quad \Big( 1 - \frac{C}{2^N \eps}
\Big) p \ge 1 \, , \quad p \le 1 + \frac{\eps}{\si^2 M} \, , \quad p \le 2 \,
.
\]
Then for $ s,r \in (0,M) $ (and $ n \ge N $, $ n |\la| \le \de $) we have
\begin{multline*}
\frac1p f_n \Big( \frac{\la}p, s \Big) + \frac1p f_n \Big( \frac{\la}p, r
 \Big) - 2\eps \le f_n(\la,s+r) \le \\
\le p f_n(p\la,s) + p f_n(p\la,r) + 2\eps \, .
\end{multline*}
Assuming in addition that $ p n |\la| \le \de $ we get
\begin{multline*}
\frac1p f_n \Big( \frac{\la}p, s+r \Big) - f_n(\la,s) - \frac{2\eps}p \le
 f_n(\la,r) \le \\
\le p f_n(p\la,s+r) - f_n(\la,s) + 2p\eps \, .
\end{multline*}
If $ s,s+r \in A $ then $ | f_n(p\la,s+r) | \le \eps $, $ | f_n \( \frac{\la}p,
s+r \) | \le \eps $, $ | f_n(\la,s) | \le \eps $, therefore $ | f_n(\la,r) |
\le 7\eps $. It remains to note that for every $ r \in (0,M) $ there exists $
s \in A $ such that $ s+r \in A $.
\end{proof}

\begin{proposition}\label{propositionC12}
Let $ X $ be splittable. Then there exists $ \si \in [0,\infty) $ such that
\[
\lim_{r\to\infty,\la\log r\to0} \frac1{ r\la^2 } \ln \Ex \exp \la \int_0^r
X_t \, \D t = \frac12 \si^2 \, .
\]
\end{proposition}

\begin{proof}
By lemma \ref{lemmaC11}, for every $ \eps $ there exist $ N $ and $ \de $ such
that
\[
\bigg| \frac1{ 2^n \la^2 } \ln \Ex \exp \la \int_0^{2^n r} X_t \, \D t -
\frac12 \si^2 r \bigg| \le \eps
\]
whenever $ n \ge N $, $ n |\la| \le \de $, $ 1 \le r \le 2 $. Dividing by $ r $
and letting $ 2^n r = s $ we get
\[
\bigg| \frac1{ s \la^2 } \ln \Ex \exp \la \int_0^s X_t \, \D t - \frac12 \si^2
\bigg| \le \eps
\]
whenever $ s \ge 2^N $ and $ |\la| \log_2 s \le \de $.
\end{proof}

\begin{lemma}\label{lemmaC13}
Let $ X $ be splittable, $ u,v,a,b \in [0,\infty) $, and $ f : (-u,v) \to
[-1,1] $ a measurable function. Then existence of two limits
\begin{align*}
\lim_{r\to\infty,\la\log r\to0} \frac1{ r\la^2 } \ln \Ex \exp \la \int_{-ur}^0
 f \Big( \frac t r \Big) X_t \, \D t = a \, , \\
\lim_{r\to\infty,\la\log r\to0} \frac1{ r\la^2 } \ln \Ex \exp \la \int_0^{vr}
 f \Big( \frac t r \Big) X_t \, \D t = b
\end{align*}
implies existence of the following limit, and the equality
\[
\lim_{r\to\infty,\la\log r\to0} \frac1{ r\la^2 } \ln \Ex \exp \la
\int_{-ur}^{vr} f \Big( \frac t r \Big) X_t \, \D t = a+b \, .
\]
\end{lemma}

\begin{proof}
Applying Lemma \ref{lemmaC4}(c) to $ ur, vr $ (and $ f(\cdot/r) $), dividing
\eqref{C6a}, \eqref{C6b} by $ r \la^2 $ and taking the limit we get (similarly
to the proof of Lemma \ref{lemmaC10})
\[
\frac1p (a+b) \le \liminf_{r\to\infty,\la\log r\to0} \frac1{ r\la^2 } \ln \Ex
\exp \la \int_{-ur}^{vr} f \Big( \frac t r \Big) X_t \, \D t \le \limsup
(\dots) \le p (a+b)
\]
for all $ p>1 $.
\end{proof}

\begin{lemma}\label{lemmaC12}
Let $ X $ be splittable. Then there exists $ \si \in [0,\infty) $ such that
\[
\lim_{r\to\infty,\la\log r\to0} \frac1{ r\la^2 } \ln \Ex \exp \la
\int_{-\infty}^\infty f \Big( \frac t r \Big) X_t \, \D t = \frac{\si^2}2 \| f
\|^2_{L_2(\R)}
\]
for every step function $ f : \R \to \R $ (that is, a \emph{finite} linear
combination of indicators of \emph{bounded} intervals).
\end{lemma}

\begin{proof}
Induction in the number $ n $ of `steps' of $ f $. For $ n=1 $ (a single
interval) it follows from Prop.~\ref{propositionC12}. The transition to $ n+1
$ is ensured by the `concatenation' Lemma \ref{lemmaC13}.
\end{proof}

\begin{theorem}\label{lemmaC133}
For every splittable process $ X $ there exist $ R \in (1,\infty) $, $ \de > 0
$ and $ M < \infty $ such that
\[
\frac1{r\la^2} \ln \Ex \exp \la \int_0^r f(t) X_t \, \D t \le M
\]
for all $ r \in [R,\infty) $, $ \la \in [-\frac\de{\ln r},0) \cup
(0,\frac\de{\ln r} ] $ and all measurable functions $ f : \R \to [-1,1] $.
\end{theorem}

\begin{proof}
Let $ \la > 0 $ (otherwise we turn to $ (-f) $).
For every $ s \in (0,1] $ Lemma \ref{easy} gives $ \Ex \exp | \int_0^s f(t)
X_t \, \D t | \le \Ex \exp |X_0| < \infty $. Lemma \ref{lemmaC2} gives a
constant $ C_1 $ such that $ \Ex \exp \la \int_0^s f(t) X_t \, \D t \le \exp
(C_1 \la^2) $ for $ |\la| \le 1 $. Lemma \ref{lemmaC4} gives us another
constant, denote it by $ C_2 $. (Note that $ C_1, C_2 $ do not depend on $ f
$.)

Defining functions $ \al_n : (0,\infty) \to [0,\infty] $ by
\begin{align*}
\al_0(\la) = \begin{cases}
 C_1 \la^2 &\text{for } \la \in (0,1], \\
 \infty &\text{for } \la \in (1,\infty),
\end{cases} \\
\al_{n+1} = (\al_n)_+ [C_2] \quad \text{for } n = 0,1,\dots
\end{align*}
we have for $ s \in (0,1] $
\[
\ln \Ex \exp \la \int_0^{2^n s} f(t) X_t \, \D t \le \al_n (\la) \quad
\text{for } \la \in (0,\infty)
\]
for $ n = 0 $, and by induction, for all $ n $; the transition from $ n $ to $
n+1 $ is based on \eqref{C6a}.

By Prop.\ \ref{lemma4},
\[
\limsup_{n\to\infty,\la n\to 0+} \frac1{2^n \la^2} \al_n(\la) < \infty \, ,
\]
which gives us $ N, \de $ and $ C $ such that $ \frac1{ 2^n \la^2 } \al_n(\la)
\le C $ whenever $ n \ge N $, $ \la n \le \de $. We have
\[
\frac1{ 2^n s \la^2 } \ln \Ex \exp \la \int_0^{2^n s} f(t) X_t \, \D t \le
\frac C s \le 2C
\]
for $ s \in [\frac12,1] $. Therefore
\[
\frac1{ r \la^2 } \ln \Ex \exp \la \int_0^r f(t) X_t \, \D t \le 2C
\]
whenever $ r \ge 2^{N-1} $ and $ \la \log_2 (2r) \le \de $.
\end{proof}

\begin{lemma}\label{lemmaC14}
Let $ X $ be splittable, and $ \si $ as in Lemma \ref{lemmaC12}. Then the set
of all $ f \in L_\infty(0,1) $ such that
\[
\lim_{r\to\infty,\la\log r\to0} \frac1{ r\la^2 } \ln \Ex \exp \la
\int_0^r f \Big( \frac t r \Big) X_t \, \D t = \frac{\si^2}2 \| f
\|^2_{L_2(0,1)}
\]
is closed.
\end{lemma}

\begin{proof}
Let $ g $ belong to this set, $ f \in L_\infty(0,1) $ and $ \| f-g \| \le \eps
< 1 $. We consider random variables
\[
X = \la \int_0^r g \Big( \frac t r \Big) X_t \, \D t \, , \quad
Y = \la \int_0^r (f-g) \Big( \frac t r \Big) X_t \, \D t \, .
\]
Theorem \ref{lemmaC133} gives us $ M $ such that
\[
\frac1{ r\la^2 } \ln \Ex \exp \la Y = \frac{ \eps^2 }{ r (\eps\la)^2 } \ln \Ex
\exp \eps\la \cdot \frac Y \eps \le M \eps^2
\]
provided that $ r $ is large enough and $ |\la| \log r $ is small enough. We
apply Lemma \ref{lemmaA0} to $ X,Y $, divide by $ r \la^2 $ and take the
limit:
\begin{multline*}
p \frac{\si^2}2 \Big\| \frac1p g \Big\|^2 - (p-1) M \Big( \frac\eps{p-1}
 \Big)^2 \le \\
\le \liminf_{r\to\infty,\la\log r\to0} \frac1{r\la^2} \ln \Ex \exp \la
 \int_0^r f \Big( \frac t r \Big) X_t \, \D t \le \limsup(\dots) \le \\
\le \frac1p \frac{\si^2}2 \| p g \|^2 + \frac{p-1}p M \Big( \frac{p}{p-1} \eps
 \Big)^2
\end{multline*}
for all $ p \in (1,\infty) $. We choose $ p = \frac1{1-\eps} $, then
\[
(1-\eps) \frac{\si^2}2 \| g \|^2 - \eps (1-\eps) M \le \liminf(\dots) \le
\limsup(\dots) \le \frac1{1-\eps} \frac{\si^2}2 \| g \|^2 + \eps M \, .
\]
Also, $ \| f-g \|_{L_2(0,1)} \le \eps $. Having such $ g $ for every $ \eps $
we get $ \liminf(\dots) = \limsup(\dots) = \frac{\si^2}2 \| f \|^2 $.
\end{proof}

\begin{proof}[Proof of Theorem \textup{\ref{theorem1}.}]
If $ f $ vanishes outside $ (0,1) $, the claim follows immediately from Lemmas
\ref{lemmaC12}, \ref{lemmaC14}. In general we take $ C $ such that $ f $
vanishes outside $ (-C,C) $ and define $ g $ by $ g(t) = f(-C+2Ct) $, then $ g
$ vanishes outside $ (0,1) $, and
\begin{multline*}
\frac1{ r\la^2 } \ln \Ex \exp \la \int f \Big( \frac t r \Big) X_t \, \D t =
 \\
= 2C \cdot \frac1{ 2Cr\la^2 } \ln \Ex \exp \la \int g \Big( \frac{t}{2Cr}
 \Big) X_t \, \D t \to 2C \cdot \frac{\si^2}2 \| g \|^2_{L_2(0,1)} =
 \frac{\si^2}2 \| f \|^2_{L_2(\R)} \, .
\end{multline*}
\end{proof}

\begin{proof}[Proof of Corollary \textup{\ref{corollary3}.}]
Let $ r_n \to \infty $, $ c_n \to \infty $, $ (c_n \log r_n )^2 / r_n \to 0 $;
we have to prove that
\[
\frac1{c_n^2} \ln \PR{ \int f \Big( \frac t{r_n} \Big) X_t \, \D t \ge c_n \si
  \| f \| \sqrt{r_n} } \to -\frac12 \quad \text{as } n \to \infty \, .
\]
By Theorem \ref{theorem1},
\[
\frac1{c_n^2} \ln \Ex \exp \frac{ \la c_n }{ \sqrt{r_n} } \int f \Big( \frac
t{r_n} \Big) X_t \, \D t \to \frac{\si^2}2 \| f \|^2 \la^2 \quad \text{as } n
\to \infty
\]
for all $ \la \in \R $. By the G\"artner-Ellis theorem \cite{Ell} (with the
scale $ c_n $ and speed $ c_n^2 $), random variables $ \frac1{ \sqrt{r_n} c_n
} \int f \( \frac t{r_n} \) X_t \, \D t $ satisfy MDP with the rate function $
x \mapsto \frac{ x^2 }{ 2\si^2 \| f \|^2 } $.
\end{proof}

\begin{proof}[Proof of Corollary \textup{\ref{corollary4}.}]
For every $ \la \in (-\infty,0) \cup (0,\infty) $, by Th.\ \ref{theorem1},
\[
\frac1{ r (\la/\sqrt r)^2 } \ln \Ex \exp \frac{\la}{\sqrt r} \int f \Big(
\frac t r \Big) X_t \, \D t \to \frac{\si^2}2 \| f \|^2 \quad \text{as } t \to
\infty \, ,
\]
that is,
\[
\Ex \exp \la \cdot \frac1{\sqrt r} \int f \Big( \frac t r \Big) X_t \, \D t
\to \exp \Big( \frac12 \si^2 \| f \|^2 \la^2 \Big) \quad \text{as } t \to
\infty \, .
\]
The weak convergence of distributions follows, see for example \cite[Sect.~30,
p.~390]{Bi}.
\end{proof}

%% file: sect3.tex
The main result of this section, Theorem \ref{th3.1}, is instrumental in
checking splittability of processes of the form $ X_t = \log |\xi_t| $ where $
\xi $ is a complex-valued Gaussian process. It is used in Sect.\
\ref{sect4}. Proposition \ref{B4} is used in Sect.\ \ref{sect6}.

Here is some background. Real-valued random processes on $ \R $ ($ X : \R
\times \Om \to \R $) are defined in the beginning of Sect.\ \ref{sect2}.
Complex-valued random processes on $ \R $ ($ X : \R \times \Om \to \C $) are
defined similarly. Let $ X $ be a complex-valued process having second moments
(that is, $ \Ex |X_t|^2 < \infty $ for all $ t \in \R $) and centered (that
is, $ \Ex X_t = 0 $ for all $ t \in \R $). Such a process is called Gaussian,
if for all $ n \in \{1,2,\dots\} $, $ t_1,\dots,t_n \in \R $ and $
a_1,\dots,a_n \in \C $ such that $ \Ex | a_1 X_{t_1} + \dots + a_n X_{t_n} |^2
= 1 $, the distribution of the random variable $ a_1 X_{t_1} + \dots + a_n
X_{t_n} $ is the standard complex normal distribution $ \ga^1 $, that is, has
the density $ z \mapsto \pi^{-1} \E^{-|z|^2} $ w.r.t.\ Lebesgue measure on $
\C $. The same holds for processes on $ (0,\infty) $, on $ \R^2 $, and in fact
on any measure space with a finite or $\si$-finite measure. Throughout this
work, all Gaussian processes are complex-valued and centered 

\begin{theorem}\label{th3.1}
Let a number $ C \in (0,\infty) $ and Gaussian processes $ \xi, \eta, \eta' $
on $ (0,\infty) $ be such that

(a) $ \xi, \eta, \eta' $ are independent;

(b) $ \eta $ and $ \eta' $ are identically distributed;

(c) for all $ u \in [0,1) $,
\[
\sum_{k=0}^\infty \( \Ex | \eta(k+u) |^2 \)^{1/2} \le C \, ;
\]
\indent
(d) for all $ u \in [0,1) $, $ n \in \{1,2,\dots\} $ and $
a_0,\dots,a_n \in \C $,
\[
\Ex \bigg| \sum_{k=0}^n a_k \( \xi(k+u) + \eta(k+u) \) \bigg|^2 \ge
\sum_{k=0}^n |a_k|^2 \, .
\]
\noindent Then
\[
\Ex \exp \int_0^\infty \ln^+ \frac{ | \xi(t) + \eta'(t) | }{ | \xi(t) +
\eta(t) | } \, \D t \le \exp \( 2\pi (C^2+C) \) \, .
\]
\end{theorem}

Here $ \ln^+ (\dots) $ means $ \max(0,\ln(\dots)) $.
The proof of Th.\ \ref{th3.1}, given at the end of this section, uses
finite-dimensional approximation (Prop.\ \ref{B4}). See also \ref{6.**}.

The standard Gaussian measure $ \ga^n $ on $ \C^n $ has the density $ x
\mapsto
\linebreak
\pi^{-n} \exp(-|x|^2) $, $ x \in \C^n $, w.r.t.\ Lebesgue measure on $
\C^n $.

Let $ \mu $ be a (finite positive Borel) measure on $ \C^n $ such that
\begin{equation}\label{2.1}
\int_{\C^n} \big| \ln |x| \big| \, \mu(\D x) < \infty \, .
\end{equation}
We denote by $ S_\mu $ the random variable
\[
S_\mu (u) = \int_{\C^n} \ln | \ip x u | \, \mu(\D x)
\]
on the probability space $ (\C^n,\ga^n) $. It is well-defined and integrable
by the Fubini theorem, since
\begin{multline*}
\iint_{ \C^n \times \C^n } \big| \ln | \ip x u | \big| \, \mu(\D x) \ga^n(\D
 u) = \\
= \int_{\C^n} \mu(\D x) \int_{\C^n} \ga^n(\D u) \, \big| \ln | \ip x u | \big|
= \int_{\C^n} \mu(\D x) \int_{\C} \ga^1(\D z) \, \big| \ln|x|+\ln|z| \big| \le
 \\
\le \int_{\C^n} \mu(\D x) \, \bigg( \big| \ln|x| \big| + \int_{\C} \ga^1(\D z)
\, \big| \ln|z| \big| \bigg) < \infty \, .
\end{multline*}

Let $ \mu,\nu $ be two measures satisfying \eqref{2.1}. We say that $ \mu $ and
$ \nu $ are projectively equal, if $ \int f \, \D\mu = \int f \, \D\nu $ for
every bounded Borel function $ f : \C^n \to \R $ such that $ f(cx) = f(x) $
for all $ x \in \C^n $, $ c \in \C \setminus \{0\} $. Clearly, every $ \mu $
is projectively equal to some $ \nu $ concentrated on the unit sphere.

\begin{sloppypar}
If $ \mu $ and $ \nu $ are projectively equal then  $ S_\mu - S_\nu $ is
nonrandom. Namely, $ S_\mu(u) - S_\nu(u) = \De_{\mu,\nu} $ for almost all
$ u $; here $ \De_{\mu,\nu} = \int_{\C^n} \ln |x| \, \mu(\D x) - \int_{\C^n}
\ln |x| \, \nu(\D x) $. Proof: the function $ f_u(x) = \ln | \ip x u | - \ln
|x| $ satisfies $ f(cx) = f(x) $. It is not bounded, however, it is both
\integrable{\mu} and \integrable{\nu} for almost every $ u $. For such $ u $
we get $ \int f_u \, \D\mu = \int f_u \, \D\nu $, thus $ S_\mu(u) -
S_\nu(u) = \De_{\mu,\nu} $.
\end{sloppypar}

If $ \mu $ and $ \nu $ are projectively equal then the corresponding centered
random variables are equal:
\[
S_\mu - \Ex S_\mu = S_\nu - \Ex S_\nu \, .
\]

Let $ \mu $ be a measure on $ \C^{m+n} $, satisfying \eqref{2.1}. Then the
random variable $ S_\mu $ is defined on the product $ (\C^{m+n},\ga^{m+n}) =
(\C^m,\ga^m) \times (\C^n,\ga^n) $ of two probability spaces. If $ \mu $ is
concentrated on $ \C^m \times \{0\} $ then $ S_\mu $ does not depend on the
second argument. If $ \mu $ is somehow close to $ \C^m \times \{0\} $, we may
expect that the second argument has a small impact on $ S_\mu $. In order to
quantify this impact we may consider the function
\[
(u,v,w) \mapsto S_\mu(u,v) - S_\mu(u,w) = \int_{\C^{m+n}} \mu ( \D y \D z
) \ln \frac{ | \ip y u + \ip z v | }{ | \ip y u + \ip z w | }
\]
as a random variable on the probability space $ (\C^m,\ga^m) \times
(\C^n,\ga^n) \times (\C^n,\ga^n) $. An upper bound for this random variable is
given in Prop.~\ref{B4} under a quite technical (but useful) condition on $
\mu $ formulated in Def.~\ref{B3}. The condition is sensitive to the choice of
a measure within a class of projectively equal measures; the conclusion is
not.

\begin{definition}\label{B3}
Let $ \mu $ be a measure on $ \C^{m+n} $, and $ A,B \in [0,\infty) $. We say
that $ \mu $ is \ABclose\ to $ \C^m \times \{0\} $, if it is an integral
convex combination (that is, $ \mu = \int_0^1 \mu_t \, \D t $) of measures $
\mu_t $ on $ \C^{m+n} $ of the following form: each $ \mu_t $ is the counting
measure on a finite set $ \{ x_1,\dots,x_K \} \subset \C^{m+n} $ (that is, $
\mu_t(A) = \# \{ k : x_k \in A \} $; both $ x_1,\dots,x_K $ and $ K $ may
depend on $ t $) satisfying the following conditions formulated in terms of $
y_k \in \C^m $, $ z_k \in \C^n $ such that $ y_k \oplus z_k = x_k $:

(a) $ |z_1| + \dots + |z_K| \le A $;

(b) $ | a_1 z_1 + \dots + a_K z_K | \le B $ whenever $ a_k \in \C $,
$ |a_k| \le 1 $ for all $ k $;

(c) $ | a_1 x_1 + \dots + a_K x_K |^2 \ge |a_1|^2 + \dots + |a_K|^2 $ for
all $ a_1,\dots,a_K \in \C $.
\end{definition}

Condition (b) makes sense for $ B < A $ (otherwise it follows from (a)).

\begin{lemma}\label{B5}
Let $ x_1,\dots,x_n \in \C^N $ satisfy $ | a_1 x_1 + \dots + a_n x_n |^2 \ge
|a_1|^2 + \dots + |a_n|^2 $ for all $ a_1,\dots,a_n \in \C $. Then
\[
\int_{\C^N} \ga^N(\D u) \prod_{k=1}^n f_k ( \ip{ x_k }{ u } ) \le
\prod_{k=1}^n \sup_{y \in \C} \int_{\C} \ga^1(\D z) \, f_k(y+z)
\]
for all measurable $ f_1,\dots,f_n : \C \to [0,\infty) $.
\end{lemma}

\begin{proof}
The image of $ \ga^N $ under the map $ \C^N \ni u \mapsto \( \ip{x_1}{u},
\dots, \ip{x_n}{u} \) \in \C^n $ is a centered Gaussian measure $ \ga $ on $
\C^n $ such that $ \int | \ip a v |^2 \, \ga(\D v) \ge |a|^2 $ for all $ a \in
\C^n $. We define another centered Gaussian measure $ \ti\ga $ on $ \C^n $ by
\[
\int | \ip a v |^2 \, \ti\ga(\D v) = \int | \ip a v |^2 \, \ga(\D v) - |a|^2
\quad \text{for all } a \in \C^n
\]
and get
\[
\int | \ip a v |^2 \, \ga(\D v) = \int | \ip a v |^2 \, \ga^n(\D v) + \int |
\ip a v |^2 \, \ti\ga(\D v) \, ,
\]
which means convolution,
\[
\ga = \ga^n * \ti\ga \, .
\]
Thus,
\begin{multline*}
\int_{\C^N} \ga^N(\D u) \prod_{k=1}^n f_k ( \ip{ x_k }{ u } ) = \int_{\C^n}
 \ga(\D v) \prod_{k=1}^n f_k (v_k) = \\
= \int_{\C^n} \ti\ga(\D v) \int_{\C^n} \ga^n(\D w) \prod_{k=1}^n f_k (v_k+w_k) =
\\
= \int_{\C^n} \ti\ga(\D v) \prod_{k=1}^n \int_{\C} \ga^1(\D z) f_k(v_k+z) \le
 \int_{\C^n} \ti\ga(\D v) \prod_{k=1}^n \sup_{y\in\C} \int_{\C} \ga^1(\D z)
 f_k(y+z) \, .
\end{multline*}
\end{proof}

\begin{proposition}\label{B4}
Let a measure $ \mu $ on $ \C^{m+n} $, satisfying \eqref{2.1}, be \ABclose\ to
$ \C^m \times \{0\} $. Then 
\begin{multline*}
\iiint_{\C^m\times\C^n\times\C^n} \!\!\!\!\! \ga^m(\D u) \ga^n(\D v) \ga^n(\D w) \, \exp
 \int_{\C^{m+n}} \mu ( \D y \D z ) \ln^+ \frac{ | \ip y u + \ip z v | }{ | \ip
 y u + \ip z w | } \le \\
\le \exp \( \pi ( A + B^2 ) \) \, .
\end{multline*}
\end{proposition}

Here $ \mu ( \D y \D z ) $ means $ \mu(\D x) $ where $ x = y \oplus z $, $ x
\in \C^{m+n} $, $ y \in \C^m $, $ z \in \C^n $.

\begin{proof}
Having $ \mu = \int_0^1 \mu_t \, \D t $ (according to Def.~\ref{B3}), we may
prove the inequality for $ \mu_t $ instead of $ \mu $ (for every $ t $), since
$ \exp \int \mu(\D y \D z) \, \ln^+(\dots) = \exp \int_0^1 \D t \int \mu_t (\D
y \D z) \, \ln^+(\dots) \le \int_0^1 \D t \exp \int \mu_t (\D y \D z) \,
\ln^+(\dots) $. The measure $ \mu_t $ being the counting measure on a finite
set of vectors $ x_k = y_k \oplus z_k $ ($ k = 1,\dots,K $, $ y_k \in \C^m $,
$ z_k \in \C^n $), we have to prove that
\begin{multline*}
\iiint_{\C^m\times\C^n\times\C^n} \ga^m(\D u) \ga^n(\D v) \ga^n(\D w) \,
 \prod_{k=1}^K \max \bigg( 1, \frac{ | \ip{y_k}u + \ip{z_k}v | }{ | \ip{y_k}u
 + \ip{z_k}w | } \bigg) \le \\
\le \exp \( \pi ( A + B^2 ) \) \, ,
\end{multline*}
given that (a) $ \sum_k |z_k| \le A $, (b) $ | \sum a_k z_k | \le B $ whenever
$ |a_k| \le 1 $, and (c) $ | \sum_k a_k y_k |^2 + | \sum_k a_k z_k |^2 \ge
\sum_k |a_k|^2 $ for all $ a_1,\dots,a_K \in \C $.

We transform the integral in $ v $ and $ w $ by the transformation
\[
(v,w) \mapsto \Big( \frac{ v-w }{ \sqrt2 }, \frac{ v+w }{ \sqrt2 } \Big)
\]
that is well-known to preserve the measure $ \ga^n \times \ga^n $ on $ \C^n
\times C^n $. Using also the evident inequality
\[
\frac{ | \ip{y_k}u + \ip{z_k}v | }{ | \ip{y_k}u + \ip{z_k}w | } \le 1 + \frac{
  | \ip{ z_k }{ v-w } | }{ | \ip{y_k}u + \ip{z_k}w | } \, ,
\]
we get
\begin{multline*}
\iint_{\C^n\times\C^n} \ga^n(\D v) \ga^n(\D w) \, \prod_{k=1}^K \max \bigg( 1,
\frac{ | \ip{y_k}u + \ip{z_k}v | }{ | \ip{y_k}u + \ip{z_k}w | } \bigg) \le \\
\le \iint_{\C^n\times\C^n} \ga^n(\D v) \ga^n(\D w) \, \prod_{k=1}^K \bigg( 1 +
 \frac{ | \ip{ z_k }{ \sqrt2 w } | }{ | \ip{y_k}u + \ip{z_k}{v+w}/\sqrt2 | }
 \bigg) = \\
= \iint_{\C^n\times\C^n} \ga^n(\D v) \ga^n(\D w) \, \prod_{k=1}^K \bigg( 1 +
 \frac{ 2 | \ip{ z_k }{ w } | }{ | \ip{ \sqrt2 y_k \oplus z_k }{ u \oplus v }
 + \ip{ z_k }{ w } | } \bigg) \, .
\end{multline*}
Applying Lemma \ref{B5} to the Gaussian measure $ \ga^m \times \ga^n $ on $
\C^{m+n} $, functions
\[
f_k (t) = 1 + \frac{ 2 | \ip{ z_k }{ w } | }{ | t + \ip{ z_k }{ w } | }
\]
and vectors $ \sqrt2 y_k \oplus z_k $ ($ k = 1,\dots,K $) we get
\begin{multline*}
\iint_{\C^m\times\C^n} \ga^m(\D u) \ga^n(\D v) \, \prod_{k=1}^K \bigg( 1 +
 \frac{ 2 | \ip{ z_k }{ w } | }{ | \ip{ \sqrt2 y_k \oplus z_k }{ u \oplus v }
 + \ip{ z_k }{ w } | } \bigg) \le \\
\le \prod_{k=1}^K \sup_{s\in\C} \int_{\C} \ga^1 (\D t) \, \Big( 1 + \frac{ 2
 | \ip{z_k}{w} | }{ | s+t | } \Big) \, .
\end{multline*}
The supremum is achieved at $ s=0 $ (which is a special case of the Anderson
inequality, see for instance \cite[Th.~1.8.5 and Cor.~1.8.6]{Bo}). Taking into
account that
\[
\int_{\C} \frac{ \ga^1(\D t) }{ |t| } = \frac{2\pi}{\pi} \int_0^\infty \frac1r
\E^{-r^2} r \D r = \sqrt\pi
\]
we get
\begin{multline*}
\iint_{\C^m\times\C^n} \ga^m(\D u) \ga^n(\D v) \, \prod_{k=1}^K \bigg( 1 +
 \frac{ 2 | \ip{ z_k }{ w } | }{ | \ip{ \sqrt2 y_k \oplus z_k }{ u \oplus v }
 + \ip{ z_k }{ w } | } \bigg) \le \\
\le \prod_{k=1}^K \( 1 + 2\sqrt\pi | \ip{ z_k }{ w } | \) \le \exp \bigg(
2\sqrt\pi \sum_{k=1}^K | \ip{ z_k }{ w } | \bigg) \, .
\end{multline*}
Thus,
\begin{multline*}
\iiint_{\C^m\times\C^n\times\C^n} \ga^m(\D u) \ga^n(\D v) \ga^n(\D w) \,
 \prod_{k=1}^K \max \bigg( 1, \frac{ | \ip{y_k}u + \ip{z_k}v | }{ | \ip{y_k}u
 + \ip{z_k}w | } \bigg) \le \\
\iiint_{\C^m\times\C^n\times\C^n} \ga^m(\D u) \ga^n(\D v) \ga^n(\D w) \,
 \prod_{k=1}^K \bigg( 1 +
 \frac{ 2 | \ip{ z_k }{ w } | }{ | \ip{ \sqrt2 y_k \oplus z_k }{ u \oplus v }
 + \ip{ z_k }{ w } | } \bigg) \\
\le \int_{\C^n} \ga^n(\D w) \, \exp \( 2\sqrt\pi f(w) \) \, ,
\end{multline*}
where $ f(w) = \sum_{k=1}^K | \ip{ z_k }{ w } | $; it is sufficient to prove
that
\[
\int_{\C^n} \ga^n(\D w) \, \exp \( 2\sqrt\pi f(w) \) \le \exp \( \pi ( A + B^2
) \) \, .
\]

It is well-known that\footnote{%
 Alternatively we could use Fernique's theorem (see for instance
 \cite[Th.~2.8.5]{Bo}), getting worse constants.}
\[
\int_{\C^n} \exp (\la f) \, \D \ga^n \le \exp \bigg( \frac14 \la^2 C^2 + \la
\int_{\C^n} f \, \D \ga^n \bigg) 
\]
for every $ f : \C^n \to \R $ satisfying the Lipschitz condition with constant
$ C $, and every $ \la \in \R $; see \cite[(1.7.8)]{Bo}. (The coefficient $
1/4 $ before $ \la^2 C^2 $ differs from the coefficient $ 1/2 $ in \cite{Bo}
because the standard Gaussian measures on $ \C^n $ and $ \R^{2n} $ have
different covariations.) It remains to check that, first, $ \int f \, \D \ga^n
\le \frac12 \sqrt\pi A $, and second, $ f $ is Lipschitz with the constant (at
most) $ B $.

The former follows from the inequality $ \sum_k |z_k| \le A $ and the fact
that
\[
\int_{\C^n} \ga^n(\D w) \, | \ip{ z_k }{ w } | = |z_k| \cdot \int_{\C} |t|
\ga^1(\D t) = |z_k| \cdot \frac{2\pi}{\pi} \int_0^\infty r \E^{-r^2} r \D r =
\frac{ \sqrt\pi }2 |z_k| \, .
\]
Finally, the gradient of $ f $ is of the form $ \sum_k a_k z_k $, $ |a_k| \le
1 $; we know that $ | \sum_k a_k z_k | \le B $.
\end{proof}

\begin{proof}[Proof of Theorem \textup{\ref{th3.1}}.]
Given $ m \in \{1,2,\dots\} $, we restrict our processes $ \xi, \eta, \eta' $
to $ (0,m) $; accordingly, we replace `$ \sum_{k=0}^\infty $' with `$
\sum_{k=0}^{m-1} $' in (c); `$ n \in \{1,2,\dots\} $' with `$ n \in
\{1,2,\dots,m-1\} $' (or just `$ n = m-1 $') in (d); and `$ \int_0^\infty
\dots \, \D t $' with `$ \int_0^m \dots \, \D t $' in the conclusion. It is
sufficient to prove this modified theorem (for all $ m $), since the limit $ m
\to \infty $ gives the original theorem.

The random variables $ \xi_t $ for $ t \in (0,m) $ span a (closed linear)
subspace $ G^{(\xi)} $ of the Hilbert space $ L_2(\Om) $ of square integrable
random variables. We choose an increasing sequence of \dimensional{n}
subspaces $ G^{(\xi)}_n \subset G^{(\xi)} $ whose union is dense in $
G^{(\xi)} $. For each $ n $ we construct a Gaussian random process $ \xi_n $
on $ (0,m) $ as follows: for every $ t \in (0,m) $ the random variable $
\xi_n(t) $ is the orthogonal projection of $ \xi(t) $ to $ G^{(\xi)}_n $. We
construct processes $ \eta_n $ in the same way. (No need to coordinate the
choice of $ G^{(\eta)}_n $ with the choice of $ G^{(\xi)}_n $.) Constructing $
\eta'_n $ we use the natural unitary correspondence between $ G^{(\eta)} $ and
$ G^{(\eta')} $ (namely, $ \eta'(t) $ corresponds to $ \eta(t) $ for each $ t
$), and construct $ \eta'_n $ such that $ \eta'_n(t) $ corresponds to $
\eta_n(t) $. Thus, for every $ n $ the three processes $ \xi_n $, $ \eta_n $,
$ \eta'_n $ satisfy Conditions (a), (b). Condition (c) is also satisfied (just
because the projection never increases norms). However, Condition (d) may be
violated.

Given $ \eps > 0 $, we'll prove the theorem with $ C $ replaced by $
(1+\eps)C $ in the conclusion (but not in Condition (c)); this is evidently
sufficient.

We define measurable sets $ A_n \subset [0,1) $ as follows: $ u \in A_n $ if
and only if the inequality
\[
(1+\eps)^2 \Ex \bigg| \sum_{k=0}^{m-1} a_k \( \xi_n(k+u) + \eta_n(k+u) \)
\bigg|^2 \ge \sum_{k=0}^{m-1} |a_k|^2
\]
holds for all $ a_0,\dots,a_{m-1} \in \C $. The left-hand side is increasing
in $ n $ (since $ \xi_n + \eta_n $ is a projection of $ \xi_{n+1} + \eta_{n+1}
$), thus $ A_1 \subset A_2 \subset \dots $ Also, $ \xi_n(k+u) + \eta_n(k+u)
\to \xi(k+u) + \eta(k+u) $ (in $ L_2 $, as $ n \to \infty $); we have an
increasing sequence of quadratic forms on $ \C^m $, and their limit is
strictly larger than $ 1 $ everywhere on the unit sphere. We see that $ u \in
A_n $ for all $ n $ large enough; that is, $ A_n \uparrow [0,1) $. We
introduce $ B_n = \bigcup_{k=0}^{m-1} (k+A_n) $ and get $ B_n \uparrow [0,m)
$.

On the other hand, $ \xi_n(t) \to \xi(t) $ (as $ n \to \infty $) a.s.\
(martingale convergence!) for each $ t $ (separately); thus, $ \xi_n(t,\om)
\to \xi(t,\om) $ for almost all pairs $ (t,\om) $. Also, $ \xi_n(t,\om)
\One_{B_n} (t) \to \xi(t,\om) $ (here $ \One_{B_n} $ is the indicator of $ B_n
$). We apply Fatou's lemma twice. First,
\[
\int_0^m \D t \, \ln^+ \frac{ | \xi(t) + \eta'(t) | }{ | \xi(t) + \eta(t) | }
\le \liminf_{n\to\infty} \int_{B_n} \D t \, \ln^+ \frac{ | \xi_n(t) +
\eta'_n(t) | }{ | \xi_n(t) + \eta_n(t) | }
\]
for almost all $ \om $; and second,
\[
\Ex \exp \int_0^m \D t \, \ln^+ \frac{ | \xi(t) + \eta'(t) | }{ | \xi(t) +
\eta(t) | } \le \liminf_{n\to\infty} \Ex \exp \int_{B_n} \D t \, \ln^+ \frac{
| \xi_n(t) + \eta'_n(t) | }{ | \xi_n(t) + \eta_n(t) | } \, .
\]
It remains to prove the inequality
\[
\Ex \exp \int_{B_n} \D t \, \ln^+ \frac{ | \xi_n(t) + \eta'_n(t) | }{ |
\xi_n(t) + \eta_n(t) | } \le \exp \( 2\pi ( (1+\eps)^2 C^2 + (1+\eps) C ) \)
\]
for all $ n $. To this end we identify both $ G^{(\xi)}_n $ and $ G^{(\eta)}_n
$ with $ \C^n $, define a measure $ \mu $ on $ G^{(\xi)}_n \oplus G^{(\eta)}_n
$ by
\[
\int \phi( x \oplus y ) \, \mu(\D x \D y) = \int_{B_n} \phi \( (1+\eps)
\xi(t), (1+\eps) \eta(t) \) \, \D t
\]
(for all bounded Borel functions $ \phi : G^{(\xi)}_n \oplus G^{(\eta)}_n \to
\R $), and apply Prop.\ \ref{B4}, taking $ A $ and $ B $ of \ref{B4} both
equal to $ (1+\eps) C $.
\end{proof}

%% file: sect4.tex
The random process considered in this section is a one-dimensional counterpart
of the random field examined in Sect.\ \ref{sect6}. Many arguments of this
section are reused in Sect.\ \ref{sect6}.

\begin{proposition}\label{propD1}
Let $ \xi = (\xi_t)_{t\in\R} $ be a stationary centered Gaussian
complex-valued random process such that
\[
\Ex \xi_s \overline{\xi_t} = \exp \( -0.5 |s-t|^2 \)
\]
for $ s,t \in \R $. Then the stationary real-valued random process $ X =
(X_t)_{t\in\R} $ defined by
\[
X_t = \ln |\xi_t| + 0.5 C_{\text{Euler}}
\]
is splittable. (Here $ C_{\text{Euler}} = 0.577\dots $ is the Euler constant.)
\end{proposition}

\begin{proof}
First,
\begin{multline*}
\Ex \exp |X_0| \le \exp ( 0.5 C_{\text{Euler}} ) \Ex \max ( |\xi_0|, 1/|\xi_0|
) = \\
= \exp ( 0.5 C_{\text{Euler}} ) \int \max(|z|,1/|z|) \ga^1(\D z) = \\
= \exp ( 0.5 C_{\text{Euler}} ) \cdot \frac{2\pi}{\pi} \int_0^\infty
\max(r,1/r) \E^{-r^2} r \, \D r < \infty
\end{multline*}
and
\begin{multline*}
\Ex X_0 = 0.5 C_{\text{Euler}} + \int \ln |z| \, \ga^1(\D z) = 0.5
 C_{\text{Euler}} + \frac{2\pi}{\pi} \int_0^\infty \ln r \cdot \E^{-r^2} r \,
 \D r = \\
= 0.5 C_{\text{Euler}} + 0.5 \int_0^\infty \ln x \cdot \E^{-x} \, \D x
= 0 \, .
\end{multline*}
We have to find $ X^0, X^-, X^+ $ satisfying Conditions (a), (b), (c) of Def.\
\ref{defI1}.

Defining a map $ \Xi : \R \to L_2(\R) $ by
\[
\Xi(t)(s) = (2/\pi)^{1/4} \exp \( -(s-t)^2 \) \, ,
\]
we get
\[
\Ex \xi_s \overline{\xi_t} = \ip{ \Xi(s) }{ \Xi(t) } \, ,
\]
since
\[
\ip{ \Xi(-t) }{ \Xi(t) } = (2/\pi)^{1/2} \int \exp \( -(s-t)^2 - (s+t)^2 \)
\D s = \exp (-(2t)^2/2) \, .
\]
We split $ \Xi $ into $ \Xi_- $ and $ \Xi_+ $,
\[
\Xi_- (t) = \Xi(t) \cdot \One_{(-\infty,0)} \, , \quad \Xi_+ (t) = \Xi(t)
\cdot \One_{(0,\infty)} \, ;
\]
as before, $ \One_A $ is the indicator of $ A $. Clearly,
\[
\Ex \xi_s \overline{\xi_t} = \ip{ \Xi_-(s) }{ \Xi_-(t) } + \ip{ \Xi_+(s) }{
\Xi_+(t) } \, .
\]
(In fact, $ \ip{ \Xi_-(s) }{ \Xi_-(t) } = \Phi(-s-t) \exp \( -0.5(s-t)^2 \) $
and $ \ip{ \Xi_+(s) }{ \Xi_+(t) } = \Phi(s+t) \exp \( -0.5(s-t)^2 \) $, where
$ \Phi(t) = (2\pi)^{-1/2} \int_{-\infty}^t \E^{-u^2/2} \, \D u $.)
We construct (on some probability space) two independent centered Gaussian
complex\nobreakdash-valued random processes $ \xi_-, \xi_+ $ such that
\begin{align*}
\Ex \xi_-(s) \overline{\xi_-(t)} &= \ip{ \Xi_-(s) }{ \Xi_-(t) } \, , \\
\Ex \xi_+(s) \overline{\xi_+(t)} &= \ip{ \Xi_+(s) }{ \Xi_+(t) } \, ,
\end{align*}
then the process $ \xi_- + \xi_+ $ is distributed like $ \xi $. Further, we
construct (on some probability space) four independent processes $ \xi_-^-,
\xi_-^+, \xi_+^-, \xi_+^+ $ such that $ \xi_-^-, \xi_-^+ $ are distributed
like $ \xi_- $ each, while $ \xi_+^-, \xi_+^+ $ --- like $ \xi_+ $.
The four processes $ \xi_-^- + \xi_+^- $, $ \xi_-^- + \xi_+^+ $, $ \xi_-^+ +
\xi_+^- $, $ \xi_-^+ + \xi_+^+ $ are distributed like $ \xi $ each. Also, the
two processes $ \xi_-^- + \xi_+^- $, $ \xi_-^+ + \xi_+^+ $ are independent. We
define $ X^0, X^-, X^+ $ by
\begin{align*}
X_t^0 &= \ln | \xi_-^-(t) + \xi_+^+(t) | + 0.5 C_{\text{Euler}} \, , \\
X_t^- &= \ln | \xi_-^-(t) + \xi_+^-(t) | + 0.5 C_{\text{Euler}} \, , \\
X_t^+ &= \ln | \xi_-^+(t) + \xi_+^+(t) | + 0.5 C_{\text{Euler}}
\end{align*}
and observe that they satisfy Conditions (a), (b); Condition (c) has to be
verified.

Using the H\"older inequality and evident symmetries,
\begin{multline*}
\!\! \Ex \exp \bigg( \int_{-\infty}^0 |X^-_t - X^0_t| \, \D t + \int_0^\infty
 |X^+_t - X^0_t| \, \D t \bigg) \le \Ex \exp \bigg( 2 \int_0^\infty
 |X^+_t - X^0_t| \, \D t \bigg) \\
= \Ex \exp \bigg( 2 \int_0^\infty \( X^+_t - X^0_t \)^- \, \D t + 2
 \int_0^\infty \( X^+_t - X^0_t \)^+ \, \D t \bigg) \le \\
\le \Ex \exp \, 4 \int_0^\infty \( X^+_t - X^0_t \)^+ \, \D t =
 \Ex \exp \, 4 \int_0^\infty \ln^+ \frac{ | \xi_-^+(t) + \xi_+^+(t) | }{ |
 \xi_-^-(t) + \xi_+^+(t) | } \, \D t = \\
= \Ex \exp \int_0^\infty \ln^+ \frac{ | \ti\xi(t) + \eta'(t) | }{ | \ti\xi(t)
 + \eta(t) | } \, \D t \, ,
\end{multline*}
where $ \ti\xi $, $ \eta $ and $ \eta' $ are defined by
\[
\ti\xi(t) = R \xi_+^+(0.25t) \, , \quad \eta(t) = R \xi_-^-(0.25t) \, ,
\quad \eta'(t) = R \xi_-^+(0.25t) \, ;
\]
the constant $ R \in (0,\infty) $ will be chosen later. Finiteness of this
expectation is ensured by Theorem \ref{th3.1} provided that Conditions
\ref{th3.1}(a,b,c,d) are satisfied by $ \ti\xi $, $ \eta $ and $ \eta' $ (for
some $ R $ and $ C $).

Conditions \ref{th3.1}(a,b) are satisfied evidently. Condition (c) is
satisfied for $ C = R \sum_{k=0}^\infty \| \Xi_-(0.25k) \| < \infty $, since
\begin{multline*}
\sum_{k=0}^\infty \( \Ex | \eta(k+u) |^2 \)^{1/2} =
 R \sum_{k=0}^\infty \( \Ex | \xi_-(0.25(k+u)) |^2 \)^{1/2} = \\
= R \sum_{k=0}^\infty \| \Xi_-(0.25(k+u)) \| \le C \, .
\end{multline*}
It remains to verify \ref{th3.1}(d). We have
\begin{multline*}
\Ex \bigg| \sum_{k=0}^n a_k \( \ti\xi(k+u) + \eta(k+u) \) \bigg|^2 = \\
= R^2 \Ex \bigg| \sum_{k=0}^n a_k \xi \( 0.25 (k+u) \) \bigg|^2 =
 R^2 \bigg\| \sum_{k=0}^n a_k \Xi (0.25 k) \bigg\|^2 \, .
\end{multline*}
Now we use unitarity of Fourier transform; taking into account that
\[
\frac1{ \sqrt{2\pi} } \int \Xi(t)(s) \E^{\I \la s} \, \D s = (2\pi)^{-1/4}
\E^{\I \la t} \E^{-\la^2/4}
\]
we get
\begin{multline*}
\Big\| \sum_k a_k \Xi(0.25k) \Big\|^2 =
\int_{-\infty}^{+\infty} \Big| \sum_k a_k \Xi(0.25k)(s) \Big|^2 \, \D s = \\
= \int_{-\infty}^{+\infty} \bigg| \sum_k a_k (2\pi)^{-1/4} \exp ( 0.25 \I
\la k ) \E^{-\la^2/4} \bigg|^2 \, \D \la = \\
= (2\pi)^{-1/2} \int_0^{8\pi} \bigg| \sum_k a_k \E^{\I \la k/4} \bigg|^2 \Big(
 \sum_{l\in\Z} \exp ( -0.5 (\la+8l\pi)^2 ) \Big) \, \D \la \ge \\
\ge \frac1{\sqrt{2\pi}} \cdot 8\pi \Big( \sum_k |a_k|^2 \Big)
\inf_{\la\in(0,8\pi)} \sum_{l\in\Z} \exp ( -0.5 (\la+8l\pi)^2 ) \, .
\end{multline*}
It remains to choose $ R $ such that
\[
R^2 \cdot \frac1{\sqrt{2\pi}} \cdot 8\pi \inf_{\la\in(0,8\pi)} \sum_{l\in\Z}
\exp ( -0.5 (\la+8l\pi)^2 ) \ge 1 \, .
\]
\end{proof}

%% file: sect5.tex
Random fields on the plane are examined in this section. The main results are
two-dimensional counterparts of Theorem \ref{theorem1} and Corollaries
\ref{corollary3}, \ref{corollary4}, formulated below
(\ref{theoremE}--\ref{5.4}) after a two-dimensional counterpart of Definition
\ref{defI1}.

Random processes on $ \R $ ($ X : \R \times \Om \to \R $) are defined in the
beginning of Sect.\ \ref{sect2}. Random fields on $ \R^2 $ ($ X : \R^2 \times
\Om \to \R $) are defined similarly. A random field $ X $ is called
stationary,\footnote{%
 In other words, homogeneous.}
if for every $ s \in \R^2 $ the shifted field $ (t,\om) \mapsto X(s+t,\om) $
is distributed like $ X $.

\begin{definition}\label{defE3}
A stationary random field $ X = (X_t)_{t\in\R^2} $ is \emph{splittable,}
if it is \splittable{C} (see below) for some $ C \in (0,\infty) $.

A stationary random field $ X = (X_t)_{t\in\R^2} $ is \emph{\splittable{C},}
if $ \Ex \exp |X_0| \le C $, $ \Ex X_0 = 0 $, and there exists (on some
probability space) a family of $ 9 $ random fields $ X^{k,l} $, $ k \in
\{-1,0,1\} $, $ l \in \{-1,0,1\} $, such that

(a1) the two triples $ (X^{-,-}, X^{-,0}, X^{-,+}) $ and $ (X^{+,-}, X^{+,0},
X^{+,+}) $ are independent,

(a2) the two triples $ (X^{-,-}, X^{0,-}, X^{+,-}) $ and $ (X^{-,+}, X^{0,+},
X^{+,+}) $ are independent;

(b1) each $ X^{k,l} $ is distributed like $ X $,

(b2) the joint distribution of the three processes $ X^{k,-}, X^{k,0}, X^{k,+}
$ does not depend on $ k \in \{-1,0,1\} $,

(b3) the joint distribution of the three processes $ X^{-,l}, X^{0,l}, X^{+,l}
$ does not depend on $ l \in \{-1,0,1\} $;

(c1) $ \Ex \exp \int_{\R} | X^{0,0}_{t_1,0} - X^{\sgn t_1,0}_{t_1,0} | \,
\D t_1 \le C $,

(c2) $ \Ex \exp \( \int_{\R} | X^{0,0}_{0,t_2} - X^{0,\sgn t_2}_{0,t_2} | \,
\D t_2 \le C $,

(c3) $ \Ex \exp \iint_{\R^2} | X^{0,0}_{t_1,t_2} - X^{\sgn t_1,0}_{t_1,t_2} -
X^{0,\sgn t_2}_{t_1,t_2} + X^{\sgn t_1,\sgn t_2}_{t_1,t_2} | \, \D t_1 \D t_2
\le C $.
\end{definition}

Of course, $ X^{-,-} $ is an abbreviation of $ X^{-1,-1} $; also, $ \sgn t $
is $ -1 $ for $ t<0 $, $ 0 $ for $ x=0 $, and $ +1 $ for $ x>0 $.

According to Def.\ \ref{defE3}, all splittable random fields are stationary.

\begin{theorem}\label{theoremE}
For every splittable random field $ X $ there exists $ \si \in [0,\infty) $
such that for every compactly supported continuous function $ f : \R^2 \to \R
$,
\[
\lim_\myatop{ r\to\infty }{ \la\log^2 r\to0 } \frac1{ r^2 \la^2 } \ln \Ex \exp
\la \iint_{\R^2} f \Big( \frac{t_1}r,\frac{t_2}r \Big) X_{t_1,t_2} \, \D t_1
\D t_2 = \frac{\si^2}2 \| f \|^2_{L_2(\R^2)} \, .
\]
\end{theorem}

For the proof see the end of this section.

\begin{corollary}\label{5.3}
Let $ X $, $ \si $ and $ f $ be as in Theorem \ref{theoremE}, and $ \si \ne 0
$. Then
\[
\lim_\myatop{ r\to\infty, c\to\infty }{ (c\log^2 r)/r \to0 } \frac1{c^2} \ln
\PR{ \iint f \Big( \frac{t_1}r,\frac{t_2}r \Big) X_{t_1,t_2} \, \D t_1 \D t_2
\ge c\si \| f \| r } = -\frac12 \, .
\]
\end{corollary}

\begin{corollary}\label{5.4}
Let $ X $, $ \si $ and $ f $ be as in Theorem \ref{theoremE}. Then the
distribution of $ r^{-1} \iint f \( \frac{t_1}r,\frac{t_2}r \) X_{t_1,t_2} \,
\D t_1 \D t_2 $ converges (as $ r \to \infty $) to the normal distribution $
N(0,\si^2 \| f \|^2) $.
\end{corollary}

Random processes are instrumental in examining random fields. We generalize
Definition \ref{defI1} and Theorem \ref{lemmaC133} to nonstationary processes
as follows.

\begin{definition}\label{defE6}
A random process $ X = (X_t)_{t\in\R} $ is \emph{splittable,} if it is
\splittable{C} (see below) for some $ C \in (0,\infty) $.

A random process $ X = (X_t)_{t\in\R} $ is \emph{\splittable{C},} if for every
$ t \in \R $, first, $ \Ex \exp |X_t| \le C $, second, $ \Ex X_t = 0 $, and
third, there exists (on some probability space) a triple of random processes $
X^0, X^-, X^+ $ such that

(a) the two processes $ X^-, X^+ $ are independent;

(b) the four processes $ X, X^0, X^-, X^+ $ are identically distributed;

(c) $ \Ex \exp \( \int_{-\infty}^t |X^-_s - X^0_s| \, \D s + \int_t^\infty
|X^+_s - X^0_s| \, \D s \) \le C $.
\end{definition}

\begin{theorem}\label{E5}
For every $ C \in (0,\infty) $ there exist $ R \in (1,\infty) $, $ \de > 0
$ and $ M < \infty $ such that for every \splittable{C} process $ X $,
\[
\frac1{r\la^2} \ln \Ex \exp \la \int_0^r f(t) X_t \, \D t \le M
\]
for all $ r \in [R,\infty) $, $ \la \in [ -\frac\de{\ln r},0) \cup
(0,\frac\de{\ln r} ] $ and all measurable functions $ f : \R \to [-1,1] $.
\end{theorem}

\begin{proof}
The proof of Theorem \ref{lemmaC133} (and Lemma \ref{lemmaC4}) needs only a
trivial adaptation.
\end{proof}

Stationary random fields can lead to nonstationary random processes as
follows. For relevance of the process $ V $ see \ref{5.88}.

\begin{lemma}\label{E6}
Let $ X $ be a splittable random field, $ X^{k,l} $ as in Def.\ \ref{defE3},
and $ f : \R^2 \to [-1,1] $ a measurable function. Then the random process $ V
$ defined by
\[
V_{t_2} = \int_\R f(t_1,t_2) \( X^{0,0}_{t_1,t_2} - X^{\sgn t_1,0}_{t_1,t_2}
\) \, \D t_1 \quad \text{for } t_2 \in \R
\]
is splittable. Moreover, if $ X $ is \splittable{C} then $ V $ is
\splittable{C} (for the same $ C $).
\end{lemma}

\begin{proof}
Let $ X $ be \splittable{C}, and $ t \in \R $. We have
\[
\Ex \exp |V_t| \le \Ex \exp \int_\R | X^{0,0}_{t_1,t} - X^{\sgn t_1,0}_{t_1,t}
| \, \D t_1 = \Ex \exp \int_{\R} | X^{0,0}_{t_1,0} - X^{\sgn t_1,0}_{t_1,0} |
\, \D t_1 \le C \, ,
\]
since the distribution of $ \( X^{0,0}_{t_1,t_2-s} - X^{\sgn
t_1,0}_{t_1,t_2-s} \)_{t_1,t_2\in\R} $ does not depend on $ s $. Clearly, $
\Ex V_t = 0 $. We define
\[
V^k_{t_2} = \int_\R f(t_1,t_2) \( X^{0,k}_{t_1,t_2-t} - X^{\sgn
t_1,k}_{t_1,t_2-t} \) \, \D t_1
\]
for $ t_2 \in \R $ and $ k \in \{-1,0,+1\} $, and check Conditions
\ref{defE6}(a,b,c).

(a) $ V^-, V^+ $ are independent, since $ V^- $ involves only $ X^{-,-},
X^{0,-}, X^{+,-} $, while $ V^+ $ involves only $ X^{-,+}, X^{0,+}, X^{+,+}
$.

(b) Each $ V^k $ is distributed like $ V $, since the distribution of $ \(
X^{0,k}_{t_1,t_2-s} - X^{\sgn t_1,k}_{t_1,t_2-s} \)_{t_1,t_2\in\R} $ does not
depend on $ s $ and $ k $. 

\smallskip
(c)
\vspace{-3mm}
\begin{multline*}
\Ex \exp \bigg( \int_{-\infty}^t |V^-_s - V^0_s| \, \D s + \int_t^\infty
 |V^+_s - V^0_s| \, \D s \bigg) = \\
= \Ex \exp \int_\R | V^{\sgn(t_2-t)}_{t_2} - V^0_{t_2} | \, \D t_2
 = \Ex \exp \int_\R | V^{\sgn t_2}_{t_2+t} - V^0_{t_2+t} | \, \D t_2 = \\
= \Ex \exp \int_\R \D t_2 \Big| \int_\R \D t_1 f(t_1,t_2+t) \( X^{0,\sgn
 t_2}_{t_1,t_2} - X^{\sgn t_1,\sgn t_2}_{t_1,t_2} - X^{0,0}_{t_1,t_2} + X^{\sgn
 t_1,0}_{t_1,t_2} \Big| \le \\
\le \Ex \exp \iint_{\R^2} | X^{0,\sgn t_2}_{t_1,t_2} - X^{\sgn t_1,\sgn
 t_2}_{t_1,t_2} - X^{0,0}_{t_1,t_2} + X^{\sgn t_1,0}_{t_1,t_2} | \, \D t_1 \D
 t_2 \le C \, .
\end{multline*}
\end{proof}

Here is a two-dimensional counterpart of Lemma \ref{lemmaC4}(b).

\begin{proposition}\label{5.88}
For every $ C \in (0,\infty) $ there exist $ R \in (1,\infty) $, $ \de > 0
$ and $ M < \infty $ such that for every \splittable{C} random field $ X $ the
following holds.

Let $ G \subset \R^2 $ be a rectangle of height $ r \ge R $, split by a
vertical line in two rectangles $ G_1, G_2 $. (In other words: $ G = [a_1,a_3]
\times [b_1,b_2] $, $ G_1 = [a_1,a_2] \times [b_1,b_2] $, $ G_2 = [a_2,a_3]
\times [b_1,b_2] $, and $ b_2 - b_1 = r \ge R $.) Let $ f : \R^2 \to [-1,1] $
be a measurable function. Then there exist (om some probability space) random
variables $ Y_1 $, $ Y_2 $ and $ Z $ such that

(a) $ Y_1,Y_2 $ are independent;

(b) $ Y_k $ is distributed like $ \iint_{G_k} f(t_1,t_2) X_{t_1,t_2} \, \D t_1
\D t_2 $ (for $ k=1,2 $);

(c) $ Y_1 + Y_2 + Z $ is distributed like $ \iint_{G} f(t_1,t_2) X_{t_1,t_2}
\, \D t_1 \D t_2 $;

(d) $ \ln \Ex \exp \la Z \le M r \la^2 $ for all $ \la \in [ -\frac\de{\ln
r},\frac\de{\ln r}] $.
\end{proposition}

\begin{proof}
Without loss of generality we assume that $ a_2 = 0 $. Using the $ 9 $ random
fields $ X^{k,l} $ of Def.\ \ref{defE3} we define
\begin{align*}
Y_1 &= \iint_{G_1} f(t_1,t_2) X_{t_1,t_2}^{-,0} \D t_1 \D t_2 = \int_{a_1}^0
 \D t_1 \int_{b_1}^{b_2} \D t_2 f(t_1,t_2) X_{t_1,t_2}^{-,0} \, , \\
Y_2 &= \iint_{G_2} f(t_1,t_2) X_{t_1,t_2}^{+,0} \D t_1 \D t_2 = \int_0^{a_3}
 \D t_1 \int_{b_1}^{b_2} \D t_2 f(t_1,t_2) X_{t_1,t_2}^{+,0} \, , \\
Z &= \iint_{G} f(t_1,t_2) \( X_{t_1,t_2}^{0,0} - X_{t_1,t_2}^{\sgn t_1,0} \)
 \D t_1 \D t_2 \, .
\end{align*}
Conditions (a), (b) hold evidently. Condition (c) holds, since $ Y_1 + Y_2 + Z
= \iint_G f(t_1,t_2) X_{t_1,t_2}^{0,0} \, \D t_1 \D t_2 $.

Lemma \ref{E6} applied to the function $ (t_1,t_2) \mapsto f(t_1,t_2) \cdot
\One_{[a_1,a_3]} (t_1) $ states that the process $ V_{t_2} = \int_{a_1}^{a_3}
\D t_1 f(t_1,t_2) \( X_{t_1,t_2}^{0,0} - X_{t_1,t_2}^{\sgn t_1,0} \) $ is
\splittable{C}. We note that $ \int_{b_1}^{b_2} V_{t_2} \, \D t_2 = Z $. 
Theorem \ref{E5} gives us $ R $, $ \de $ and $ M $ (dependent on $ C $ only)
such that $ \ln \Ex \exp \la Z \le M r \la^2 $ provided that $ b_2 - b_1 = r
\ge R $ and $ \la \in [ -\frac\de{\ln r},\frac\de{\ln r}] $.
\end{proof}

\begin{remark}
The same holds for a rectangle split by a \emph{horizontal} line, since Def.\
\ref{defE3} is symmetric w.r.t.\ the interchange of the two coordinates, $
(X_{t_1,t_2})_{t_1,t_2} \mapsto (X_{t_2,t_1})_{t_1,t_2} $.
\end{remark}

In the following proposition, the simplest case $ d=0 $ amounts to Prop.\
\ref{lemma4} plus a part (upper bound) of Prop.\ \ref{prop1}. The case $ d=1 $
is used in the proof of Th.\ \ref{theoremE} in two ways: via \ref{5.8}, and
via \ref{lemma311} - \ref{5.13} - \ref{5.14} - \ref{5.15} - \ref{5.18} -
\ref{5.19}.

\begin{proposition}\label{lemma31}
Let functions $ f_n, g_n : (0,\infty) \to [0,\infty] $ satisfy
\[
f_{n+1}(\la) \le \frac2p f_n (p\la) + \frac{p-1}p g_n \Big( \frac{p}{p-1} \la
\Big)
\]
for all $ n $, $ \la $ and $ p \in (1,\infty) $, and
\[
\limsup_{\la\to0+} \frac1{\la^2} f_n(\la) < \infty
\]
for all $ n $. Let $ d \in \{0,1,2,\dots\} $. If
\[
\limsup_{n\to\infty,\la n^d\to0+} \frac1{ (2\theta)^n \la^2 } g_n(\la) < \infty
\quad \text{for some } \theta \in (0,1) \, ,
\]
then
\[
\limsup_{n\to\infty,\la n^{d+1}\to0+} \frac1{ 2^n \la^2 } f_n(\la) =
\lim_{n\to\infty} \limsup_{\la\to0+} \frac1{2^n \la^2} f_n(\la) < \infty \, .
\]
\end{proposition}

\begin{proof}
We take $ \theta \in (0,1) $, $ \de > 0 $ and $ N $ such that
\[
\sup_{n\ge N,\la n^d\le\de} \frac1{ (2\theta)^n \la^2 } g_n(\la) = M < \infty
\, .
\]
Given $ n > N $, we define functions $ h_N, h_{N+1}, \dots, h_n : (0,\infty)
\to [0,\infty] $ by
\[
h_k(\la) = f_k ( \de n^{-d} \la ) \, .
\]
Then
\begin{multline*}
h_{k+1}(\la) = f_{k+1} ( \de n^{-d} \la ) \le \frac2p f_k(p \de n^{-d} \la) +
 \frac{p-1}p g_k \Big( \frac{p}{p-1} \de n^{-d} \la \Big)  \le \\
\le \frac2p h_k(p \la) + \frac{p-1}p M (2\theta)^k \Big( \frac{p}{p-1} \de
 n^{-d} \la \Big)^2 = \\
= \frac2p h_k(p \la) + \frac{p}{p-1} M (2\theta)^k ( \de n^{-d} )^2 \la^2
\end{multline*}
for all $ p $ such that $ \frac{p}{p-1} \de n^{-d} \la \cdot n^d \le \de $,
that is, $ \frac{p}{p-1} \la \le 1 $. It means that
\[
h_{k+1} \le (h_k)_+ [C_k] \, , \quad \text{where } C_k = M (2\theta)^k ( \de
n^{-d} )^2 \quad (k=N,N+1,\dots,n-1) \, .
\]
By Prop.\ \ref{lemma44}(a),
\[
\frac1{ 2^{n-N} } h_n(\la) \le (1-\eps) h_N \Big( \frac{ \la }{ 1-\eps } \Big)
+ \frac{ \la^2 }{ 2\eps (1-\sqrt\theta)^2 } M ( \de n^{-d} )^2 (2\theta)^N
\]
for $ \eps \in (0,1) $ and $ \la \le \eps \theta^{(n-N-1)/2} (1-\sqrt\theta)
$. That is,
\[
\frac1{ 2^n } f_n(\la) \le (1-\eps) \frac1{ 2^N } f_N \Big( \frac{ \la }{
1-\eps } \Big) + \frac{ M \la^2 \theta^N }{ 2 \eps (1-\sqrt\theta)^2 }
\]
for $ \la \le \de n^{-d} \eps \theta^{(n-N-1)/2} (1-\sqrt\theta) $. (Now we may
forget the functions $ h_N, \dots, h_n $ and release $ n $.) Therefore
\begin{equation}\label{C2}
\limsup_{n\to\infty, \la n^d \theta^{-n/2} \to 0+} \frac1{ 2^n \la^2 } f_n(\la)
\le \frac1{1-\eps} \frac1{ 2^N } \limsup_{\la\to0+} \frac1{\la^2} f_N(\la) +
 \frac{ M \theta^N }{ 2 \eps (1-\sqrt\theta)^2 } \, .
\end{equation}
On the other hand, by Prop.\ \ref{lemma44}(b),
\[
\frac1{ 2^{n-N} } h_n(\la) \le (1-\eps) h_N \Big( \frac{ \mu }{ 1-\eps } \Big)
+ \bigg( \frac{ \mu^2 }{ 2\eps (1-\sqrt\theta)^2 } + \frac\la2 \frac{
\theta^{m-N} }{ 1-\theta } \bigg) M ( \de n^{-d} )^2 (2\theta)^N
\]
whenever $ m \in \{ N+1,N+2,\dots,n-1 \} $ is such that $ \la < \frac1{n-m} $
and $ \mu \le \eps \theta^{(m-N-1)/2} (1-\sqrt\theta) $, where $ \mu = \frac{
\la }{ 1 - (n-m) \la } $. That is,
\begin{multline*}
\frac1{ 2^n } f_n(\la) \le (1-\eps) \frac1{ 2^N } f_N \Big( \frac{ \la }{
 1-\eps } \cdot \frac1{ 1 - (n-m) n^d \de^{-1} \la } \Big) + \\
+ \bigg( \frac{ \la^2 }{ ( 1 - (n-m) n^d \de^{-1} \la )^2 2\eps
 (1-\sqrt\theta)^2 } + \frac{ \de \la }{ 2n^d } \frac{ \theta^{m-N} }{
 1-\theta } \bigg) M \theta^N
\end{multline*}
provided that $ \la < \de n^{-d} \frac1{n-m} $ and $ \frac{ \la }{ 1 -
(n-m) n^d \de^{-1} \la } < \de n^{-d} \eps \theta^{(m-N-1)/2} (1-\sqrt\theta)
$. Assuming $ \theta^{n-N} \ll \la n^d \ll \frac1n $, we choose $ m $ such that
\[
( \la n^d )^2 \ll \theta^{m-N} \ll \la n^d \, ,
\]
then $ (n-m) n^d \de^{-1} \la \ll 1 $, $ \de \la n^{-d} \theta^{m-N} \ll \la^2
$ and therefore
\begin{multline*}
\limsup_{n\to\infty, \la n^{d+1} \to 0+, \la n^d \theta^{-n} \to \infty}
 \frac1{ 2^n \la^2 } f_n(\la) \le \\
\le \frac1{1-\eps} \frac1{ 2^N } \limsup_{\la\to0+} \frac1{\la^2} f_N(\la) +
 \frac{ M \theta^N }{ 2 \eps (1-\sqrt\theta)^2 } \, .
\end{multline*}
In combination with \eqref{C2} it gives us
\[
\limsup_{n\to\infty, \la n^{d+1} \to 0+} \frac1{ 2^n \la^2 } f_n(\la)
\le \frac1{1-\eps} \frac1{ 2^N } \limsup_{\la\to0+} \frac1{\la^2} f_N(\la) +
 \frac{ M \theta^N }{ 2 \eps (1-\sqrt\theta)^2 } < \infty \, .
\]
This inequality holds for all $ N $ large enough; we take $ N \to \infty $ and
then $ \eps \to 0 $:
\[
\limsup_{n\to\infty, \la n^{d+1} \to 0+} \frac1{ 2^n \la^2 } f_n(\la)
\le \liminf_{N\to\infty} \limsup_{\la\to0+} \frac1{2^N \la^2} f_N(\la) \, .
\]
However, evidently
\begin{multline*}
\liminf_{N\to\infty} \limsup_{\la\to0+} \frac1{2^N \la^2} f_N(\la) \le
 \limsup_{N\to\infty} \limsup_{\la\to0+} \frac1{2^N \la^2} f_N(\la) \le \\
\le \limsup_{n\to\infty, \la n^{d+1} \to 0+} \frac1{ 2^n \la^2 } f_n(\la) \, ,
\end{multline*}
which means that all these numbers are equal.
\end{proof}

Here is the corresponding lower bound.

\begin{proposition}\label{5.11}
Let functions $ f_n, g_n : (0,\infty) \to [0,\infty] $ satisfy
\[
f_{n+1}(\la) \ge 2p f_n \Big( \frac{\la}p \Big) - (p-1) g_n \Big(
\frac{\la}{p-1} \Big)
\]
for all $ n $, $ \la $ and $ p \in (1,\infty) $. Let $ d \in \{0,1,2,\dots\}
$. If
\[
\limsup_{n\to\infty,\la n^d\to0+} \frac1{ (2\theta)^n \la^2 } g_n(\la) < \infty
\quad \text{for some } \theta \in (0,1) \, ,
\]
then
\[
\liminf_{n\to\infty,\la n^{d+1}\to0+} \frac1{ 2^n \la^2 } f_n(\la) =
\lim_{n\to\infty} \liminf_{\la\to0+} \frac1{2^n \la^2} f_n(\la) \, .
\]
\end{proposition}

\begin{proof}
Let $ \theta, \de, N, n $, $ h_N,\dots,h_n $ and $ C_N,\dots,C_{n-1} $ be as
in the proof of Prop.\ \ref{lemma31}, then $ h_{k+1} \ge (h_k)_- [C_k] $. We
use Prop.\ \ref{lemma444}(a,b) instead of Prop.\ \ref{lemma44}(a,b) and get
(similarly to \ref{lemma31})
\[
\liminf_{n\to\infty, \la n^{d+1} \to 0+} \frac1{ 2^n \la^2 } f_n(\la)
\ge \limsup_{N\to\infty} \liminf_{\la\to0+} \frac1{2^N \la^2} f_N(\la) \, .
\]
However, evidently
\begin{multline*}
\limsup_{N\to\infty} \liminf_{\la\to0+} \frac1{2^N \la^2} f_N(\la) \ge
 \liminf_{N\to\infty} \liminf_{\la\to0+} \frac1{2^N \la^2} f_N(\la) \ge \\
\ge \liminf_{n\to\infty, \la n^{d+1} \to 0+} \frac1{ 2^n \la^2 } f_n(\la) \, ,
\end{multline*}
which means that all these numbers are equal.
\end{proof}

Now we combine the lower and upper bounds, getting a generalization of Prop.\
\ref{prop1}.

\begin{proposition}\label{lemma311}
Let functions $ f_n, g_n : (0,\infty) \to [0,\infty] $ satisfy
\[
2p f_n \Big( \frac{\la}p \Big) - (p-1) g_n \Big( \frac{\la}{p-1} \Big)
\le f_{n+1}(\la) \le \frac2p f_n (p\la) + \frac{p-1}p g_n \Big( \frac{p}{p-1}
 \la \Big)
\]
for all $ n $, $ \la $ and $ p \in (1,\infty) $. Assume existence of the limit
\[
\lim_{\la\to0+} \frac1{\la^2} f_n(\la) \in [0,\infty)
\]
for all $ n $. Let $ d \in \{0,1,2,\dots\} $. If
\[
\limsup_{n\to\infty,\la n^d\to0+} \frac1{ (2\theta)^n \la^2 } g_n(\la) < \infty
\quad \text{for some } \theta \in (0,1) \, ,
\]
then the following limit exists:
\[
\lim_{n\to\infty,\la n^{d+1}\to0+} \frac1{ 2^n \la^2 } f_n(\la) \in [0,\infty)
\, .
\]
\end{proposition}

\begin{proof}
The limit exists, since $ \limsup $ and $ \liminf $ are equal; indeed, by
Prop.\ \ref{lemma31} and \ref{5.11} they are both equal to
\[
\lim_{n\to\infty} \lim_{\la\to0+} \frac1{ 2^n \la^2 } f_n(\la) < \infty \, .
\]
\end{proof}

Here are two-dimensional counterparts of \ref{propC}, \ref{lemmaC10},
\ref{lemmaC11}, \ref{propositionC12} and \ref{lemmaC12}.

\begin{proposition}\label{5.13}
For every splittable random field $ X $ and all $ r_1,r_2 \in (0,\infty) $
the following limit exists:
\[
\lim_{n\to\infty,\la n^2\to0} \frac1{ 4^n \la^2 } \ln \Ex \exp \la \int_0^{2^n
r_1} \D t_1 \int_0^{2^n r_2} \D t_2 X_{t_1,t_2} \, .
\]
\end{proposition}

\begin{proof}
Let $ s_1,s_2 \in (0,1] $. We define functions $ f_n : \R \to [0,\infty] $ by
\[
f_n(\la) = \ln \Ex \exp \la \iint_{G_n} X_{t_1,t_2} \, \D t_1 \D t_2 \, ,
\]
rectangles $ G_n $ of area $ 2^n s_1 s_2 $ being
\[
G_{2k} = [0,2^k s_1] \times [0,2^k s_2] \, , \quad
G_{2k+1} = [0,2^{k+1} s_1] \times [0,2^k s_2] \, .
\]
Note that $ f_n(\la) \ge 0 $ by Lemma \ref{lemmaA**}, and $ f_n(\la) < \infty
$ for $ |\la| \le 2^{-n} $ by (the evident two-dimensional counterpart of)
Lemma \ref{easy}.

The rectangle $ G_{n+1} $ can be split in two rectangles congruent to $ G_n $
(moreover, shifts of $ G_n $) by a line of length $ \le 2^{(n+1)/2} $. Prop.\
\ref{5.88} (in combination with Lemma \ref{lemmaA0}) provides $ R $, $ \de $
and $ M $ such that for all $ \la \in \R $ and $ p \in (1,\infty) $,
\begin{equation}\label{5.13a}
2p f_n \Big( \frac{\la}p \Big) - (p-1) g_n \Big( \frac{\la}{p-1} \Big)
\le f_{n+1}(\la) \le \frac2p f_n (p\la) + \frac{p-1}p g_n \Big( \frac{p}{p-1}
\la \Big) \, ,
\end{equation}
where $ g_n : \R \to [0,\infty] $ satisfy
\[
g_n(\la) \le M \cdot 2^{(n+1)/2} \la^2 \quad \text{whenever } 2^{(n+1)/2} \ge
R \text{ and } |\la| \frac{n+1}2 \ln 2 \le \de \, .
\]
We have
\[
\limsup_{n\to\infty,\la n\to0} \frac1{ (\sqrt2)^n \la^2 } g_n(\la) \le M \cdot
\sqrt2 < \infty \, .
\]
Prop.\ \ref{lemma311} (for $ d=1 $) ensures existence of the limit
\[
\lim_{n\to\infty,\la n^2\to0+} \frac1{ 2^n \la^2 } f_n(\la) \in [0,\infty) \,
.
\]
The same argument applied to $ (-X) $ gives us $ \lim_{n\to\infty,\la
n^2\to0-} \frac{ f_n(\la) }{ 2^n \la^2 } $. The two limits are equal, since
both are equal to $ \lim_{n\to\infty} \( 2^{-n-1} f''_n(0) \) $.

Given $ r_1,r_2 \in (0,\infty) $, we take $ m $ and $ s_1,s_2 \in (0,1] $ such
that $ r_1 = 2^m s_1 $, $ r_2 = 2^m s_2 $ and note that
\[
\frac1{ 4^n \la^2 } \ln \Ex \exp \la \int_0^{2^n r_1} \D t_1 \int_0^{2^n r_2}
\D t_2 X_{t_1,t_2} = 2^{2m} \cdot \frac1{ 2^{2(n+m)} \la^2 } f_{2(n+m)} (\la)
\, .
\]
\end{proof}

\begin{lemma}\label{5.14}
For every splittable random field $ X $ there exist $ \si \in [0,\infty) $
such that
\[
\lim_{n\to\infty,\la n^2\to0} \frac1{ 4^n \la^2 } \ln \Ex \exp \la \int_0^{2^n
r_1} \D t_1 \int_0^{2^n r_2} \D t_2 X_{t_1,t_2} = \frac12 \si^2 r_1 r_2
\]
for all $ r_1,r_2 \in (0,\infty) $.
\end{lemma}

\begin{proof}
Similarly to the proof of \ref{lemmaC10} we denote the limit (given by
\ref{5.13}) by $ \phi(r_1,r_2) $, use Prop.\ \ref{5.88} for proving additivity
in $ r_2 $ (for any given $ r_1 $), as well as in $ r_1 $ (for any given $ r_2
$), and conclude that $ \phi(r_1,r_2) = c r_1 r_2 $.
\end{proof}

\begin{lemma}\label{5.15}
The convergence in Lemma \ref{5.14} is uniform in $ r_1,r_2 \in (0,M) $ for
every $ M \in (0,\infty) $.
\end{lemma}

\begin{proof}
Similarly to the proof of \ref{lemmaC11} we use the (approximate) additivity
in $ r_2 $ for any given $ r_1 $ (ensured by Prop.\ \ref{5.88}) and get the
uniform convergence for $ (r_1,r_2) \in A \times (0,M) $ for some measurable
set $ A \subset (0,2M) $ of Lebesgue measure $ > \frac32 M $. Then we use
additivity in $ r_1 $ (for any given $ r_2 $).
\end{proof}

\begin{proposition}\label{5.18}
For every splittable random field $ X $ there exist $ \si \in [0,\infty) $
such that for every $ a \in (0,\infty) $
\[
\lim_\myatop{ r\to\infty }{ \la \log^2 r \to0 } \frac1{ ar^2 \la^2 } \ln \Ex
\exp \la \int_0^r \D t_1 \int_0^{ar} \D t_2 X_{t_1,t_2} = \frac12 \si^2 \, .
\]
\end{proposition}

\begin{proof}
Follows from Lemma \ref{5.15} similarly to the proof of \ref{propositionC12}.
\end{proof}

\begin{lemma}\label{5.19}
For every splittable random field $ X $ there exist $ \si \in [0,\infty) $
such that
\[
\lim_\myatop{ r\to\infty }{ \la \log^2 r \to0 } \frac1{ r^2 \la^2 } \ln \Ex
\exp \la \iint_{\R^2} f \Big( \frac{t_1}r,\frac{t_2}r \Big) X_{t_1,t_2} \, \D
t_1 \D t_2 = \frac{\si^2}2 \| f \|_{L_2(\R^2)}^2
\]
for every function $ f : \R^2 \to \R $ that is a (finite) linear combination
of indicators of (bounded) rectangles.
\end{lemma}

\begin{proof}
Similarly to the proof of \ref{lemmaC12}, we apply the concatenation
argument (supported by \ref{5.88}) in $ t_2 $, thus getting (via \ref{5.18})
functions $ f $ of the form $ f(t_1,t_2) = \One_{(a,b)}(t_1) f_2(t_2) $ where
$ f_2 $ runs over step functions. Then we apply the concatenation argument
again, this time in $ t_1 $.
\end{proof}

Here is a two-dimensional counterpart of Th.\ \ref{E5} (see also
\ref{lemmaC133}).

\begin{theorem}\label{5.8}
For every $ C \in (0,\infty) $ there exist $ R \in (1,\infty) $, $ \de > 0
$ and $ M < \infty $ such that for every \splittable{C} stationary random
field $ X $,
\[
\frac1{r^2 \la^2} \ln \Ex \exp \la \int_0^r \D t_1 \int_0^r \D t_2 f(t_1,t_2)
X_{t_1,t_2} \le M 
\]
for all $ r \in [R,\infty) $, $ \la \in [-\frac\de{\ln^2 r},0) \cup
(0,\frac\de{\ln^2 r} ] $ and all measurable functions $ f : \R^2 \to [-1,1]
$.
\end{theorem}

\begin{proof}
We consider only $ \la > 0 $ (otherwise turn to $ (-f) $). We define functions
$ \al_n : (0,\infty) \to [0,\infty] $ by
\[
\al_n (\la) = \sup_{X,f} \ln \Ex \exp \la \iint_{G_n} f(t_1,t_2) X_{t_1,t_2}
\, \D t_1 \D t_2 \, ;
\]
here the supremum is taken over all measurable $ f : \R^2 \to [-1,1] $ and all 
\splittable{C} $ X $, and rectangles $ G_n $ of area $ 2^n $ are as in the
proof of \ref{5.13} (for $ s_1=s_2=1 $). We have $ \al_n(2^{-n}) \le \Ex |X_0|
\le C < \infty $ by (the evident two-dimensional counterpart of) Lemma
\ref{easy}. It follows by Lemma \ref{lemmaC2} that
\[
\limsup_{\la\to0+} \frac1{\la^2} \al_n(\la) < \infty \quad \text{for all } n
\, .
\]
Using the upper bound of \eqref{5.13a} we get
\begin{gather*}
\al_{n+1}(\la) \le \frac2p \al_n (p\la) + \frac{p-1}p g_n \Big( \frac{p}{p-1}
 \la \Big) \, , \\
\limsup_{n\to\infty,\la n\to0+} \frac1{ (\sqrt2)^n \la^2 } g_n(\la) < \infty
\, .
\end{gather*}
By Prop.\ \ref{lemma31} (for $ d=1 $),
\[
\limsup_{n\to\infty,\la n^2 \to0+} \frac1{ 2^n \la^2 } \al_n(\la) < \infty \,
.
\]
Now the proof is finalized by the argument of the last paragraph of the proof
of \ref{lemmaC133}.
\end{proof}

\begin{proof}[Proof of Theorem \textup{\ref{theoremE}.}]
Similarly to the proof of Lemma \ref{lemmaC14}, using Th.\ \ref{5.8} we check
that the set of all $ f \in L_\infty \( (0,1) \times (0,1) \) $ satisfying
\[
\lim_\myatop{ r\to\infty }{ \la\log^2 r\to0 } \frac1{r^2 \la^2} \ln \Ex \exp
\la \int_0^r \D t_1 \int_0^r \D t_2 f \Big( \frac{t_1}r, \frac{t_2}r \Big)
X_{t_1,t_2} = \frac{ \si^2 }2 \| f \|^2_{L_2((0,1) \times (0,1))}
\]
is closed. Due to \ref{5.19}, it contains all continuous functions. The rest
of the proof is similar to the proof of Th.\ \ref{theorem1}.
\end{proof}

Corollaries \ref{5.3} and \ref{5.4} follow from Th.\ \ref{theoremE} in the
same way as Corollaries \ref{corollary3} and \ref{corollary4} from Th.\
\ref{theorem1}.

%% file: sect6.tex
Splittability of the random field examined in this section is significant for
Sect.\ \ref{sect7}.

We consider a centered Gaussian complex-valued random field $ \xi =
(\xi_t)_{t\in\R^2} $ such that
\[
\Ex \xi_s \overline{\xi_t} = \exp ( -\I s \wedge t - 0.5 |s-t|^2 )
\]
(here $ s \wedge t = (s_1,s_2) \wedge (t_1,t_2) = s_1 t_2 - s_2 t_1 $), and
the real-valued random field $ X = (X_t)_{t\in\R^2} $ defined by
\[
X_t = \al \ln |\xi_t| + 0.5 \al C_{\text{Euler}} \, ;
\]
an absolute constant $ \al \in (0,1] $ will be chosen later. The field
satisfies $ \Ex \exp |X_0| < \infty $ and $ \Ex X_0 = 0 $ (see the proof of
Prop.~\ref{propD1}).

\begin{lemma}
$ X $ is stationary.
\end{lemma}

\begin{proof}
For every $ r \in \R^2 $ the random field
\[
\( \exp ( \I t \wedge r ) \xi_{t+r} \)_{t\in\R^2}
\]
is distributed like $ \xi $, since
\begin{multline*}
\Ex \exp ( \I s \wedge r ) \xi_{s+r} \overline{ \exp ( \I t \wedge r )
 \xi_{t+r} } = \\
= \exp \( \I s \wedge r - \I t \wedge r - \I (s+r) \wedge (t+r) - 0.5 |s-t|^2
\) = \\
= \exp ( -\I s \wedge t - 0.5 |s-t|^2 ) = \Ex \xi_s \overline{\xi_t} \, .
\end{multline*}
Therefore the random field $ \( \al \ln |\xi_{t+r}| + 0.5 \al C_{\text{Euler}}
\)_{t\in\R^2} = (X_{t+r})_{t\in\R^2} $ is distributed like $ X $.
\end{proof}

We define a map $ \Xi : \R^2 \to L_2(\R^2) $ by
\[
\Xi(t)(s) = \frac1{\sqrt\pi} \exp \( - \I t \wedge s - 0.5 |s-t|^2 \) \, .
\]

The following fact is well-known, see for instance \cite[(7.49)]{KS}.

\begin{lemma}
$ \Ex \xi_s \overline{\xi_t} = \ip{ \Xi(s) }{ \Xi(t) } $ for $ s,t \in \R^2
$.
\end{lemma}

\begin{proof}
Using the equality
\[
\frac1\pi \int \exp ( \I u \wedge r ) \exp (-|r|^2) \, \D r = \exp ( -0.25
|u|^2 )
\]
we get
\begin{multline*}
\ip{ \Xi(s) }{ \Xi(t) } = \int \Xi(s)(r) \overline{ \Xi(t)(r) } \, \D r = \\
= \frac1\pi \int \exp ( -\I s \wedge r - 0.5 |s-r|^2 + \I t \wedge r - 0.5
 |t-r|^2 ) \, \D r = \\
= \frac1\pi \int \exp \( -\I (s-t) \wedge r - | r - 0.5(s+t) |^2 - 0.25
 |s-t|^2 \) \, \D r = \\
= \exp \( - 0.25 |s-t|^2 \) \frac1\pi \int \exp \( -\I (s-t) \wedge ( r +
 0.5(s+t) ) - |r|^2 \) \, \D r = \\
= \exp \( - 0.25 |s-t|^2 - 0.5 \I (s-t) \wedge (s+t) \) \frac1\pi \int \exp \(
 -\I (s-t) \wedge r \) \exp ( -|r|^2 ) \, \D r = \\
= \exp \( - 0.25 |s-t|^2 - \I s \wedge t \) \exp \( - 0.25 |s-t|^2
\) = \Ex \xi_s \overline{\xi_t} \, .
\end{multline*}
\end{proof}

We introduce $ \Xi_{-,-}, \Xi_{-,+}, \Xi_{+,-}, \Xi_{+,+} : \R^2 \to L_2(\R^2)
$ by
\begin{gather*}
\Xi_{a_1,a_2} (t) = \Xi(t) \cdot I_{a_1,a_2} \, , \quad \text{where} \\
I_{a_1,a_2} (t_1,t_2) = \begin{cases}
 1 &\text{if } \sgn t_1 = a_1 \text{ and } \sgn t_2 = a_2, \\
 0 &\text{otherwise}
\end{cases}
\end{gather*}
for $ a_1, a_2 = \pm 1 $. Clearly,
\[
\Ex \xi_s \overline{\xi_t} = \sum_{a_1,a_2} \ip{ \Xi_{a_1,a_2}(s) }{
\Xi_{a_1,a_2}(t) } \, .
\]
We construct (on some probability space) four independent centered Gaussian
complex-valued random fields $ \xi_{a_1,a_2} $ such that
\[
\Ex \xi_{a_1,a_2}(s) \overline{\xi_{a_1,a_2}(t)} = \ip{ \Xi_{a_1,a_2}(s) }{
\Xi_{a_1,a_2}(t) }
\]
for $ a_1, a_2 = \pm 1 $; then the field $ \sum_{a_1,a_2} \xi_{a_1,a_2} $ is
distributed like $ \xi $. Further, we construct (on some probability space)
$ 16 $ random fields $ \xi_{a_1,a_2}^{b_1,b_2} $ ($ a_1, a_2, b_1, b_2
= \pm 1 $) such that
\begin{equation}\label{6.*}
\begin{gathered}
\text{the $ 16 $ fields $ \xi_{a_1,a_2}^{b_1,b_2} $ are independent;}\\
\text{each $ \xi_{a_1,a_2}^{b_1,b_2} $ is distributed like $ \xi_{a_1,a_2} $.}
\end{gathered}
\end{equation}
We define $ 9 $ random fields $ X^{k,l} $ ($ k, l \in \{-1,0,1\} $),
distributed like $ X $ each, by
\begin{align*}
X_t^{0,0} &= \al \ln \bigg| \sum_{a_1,a_2=\pm1} \xi_{a_1,a_2}^{a_1,a_2} (t)
 \bigg| + 0.5 \al C_{\text{Euler}} \, , \\
X_t^{b_1,0} &= \al \ln \bigg| \sum_{a_1,a_2=\pm1} \xi_{a_1,a_2}^{b_1,a_2} (t)
 \bigg| + 0.5 \al C_{\text{Euler}} \, , \\
X_t^{0,b_2} &= \al \ln \bigg| \sum_{a_1,a_2=\pm1} \xi_{a_1,a_2}^{a_1,b_2} (t)
 \bigg| + 0.5 \al C_{\text{Euler}} \, , \\
X_t^{b_1,b_2} &= \al \ln \bigg| \sum_{a_1,a_2=\pm1} \xi_{a_1,a_2}^{b_1,b_2} (t)
 \bigg| + 0.5 \al C_{\text{Euler}}
\end{align*}
for $ b_1, b_2 \in \{-1,1\} $. In order to prove splittability of $ X $ we
check Conditions (a1)--(c3) of Def.\ \ref{defE3}.
Condition \ref{defE3}(a1) is satisfied (as well as (a2)), since the fields $
X^{-,-}, X^{-,0}, X^{-,+} $ involve $ \xi_{a_1,a_2}^{b_1,b_2} $ for $ b_1=-1 $
only, while $ X^{+,-}, X^{+,0}, X^{+,+} $ --- for $ b_1=+1 $ only.

\begin{lemma}
Conditions \ref{defE3}(b1,b2,b3) are satisfied.
\end{lemma}

\begin{proof}
According to \eqref{6.*}, the four fields $ \xi_{a_1,a_2}^{+,+} $ are
independent and distributed like $ \xi_{a_1,a_2} $; therefore $ X^{+,+} \sim X
$ (here `$ \sim $' means `distributed like'). The same holds for the four
fields $ \xi_{a_1,a_2}^{a_1,a_2} $, thus $ X^{0,0} \sim X $. Similarly, $
X^{k,l} \sim X $ for all $ k,l \in \{-1,0,1\} $, which is (b1).

According to \eqref{6.*}, the 8-component field $
(\xi_{a_1,a_2}^{-,b_2})_{a_1,a_2,b_2=\pm1} $ is distributed like $
(\xi_{a_1,a_2}^{+,b_2})_{a_1,a_2,b_2=\pm1} $. Therefore the triple $ (X^{-,-},
X^{-,0}, X^{-,+}) $ is distributed like $ (X^{+,-},X^{+,0}, X^{+,+}) $ (which
is a part of (b2)). Also, $ (\xi_{a_1,a_2}^{a_1,b_2})_{a_1,a_2,b_2=\pm1} \sim 
(\xi_{a_1,a_2}^{+,b_2})_{a_1,a_2,b_2=\pm1} $, since the transformation $
(a_1,a_2,b_1,b_2) \mapsto (a_1,a_2,a_1 b_1, b_2) $ leaves \eqref{6.*}
invariant. Therefore $ (X^{0,-},X^{0,0}, X^{0,+}) \sim
(X^{+,-},X^{+,0},X^{+,+}) $, which completes the proof of (b2); (b3) is
similar.
\end{proof}

\begin{lemma}
Conditions \ref{defE3}(c1,c2) are satisfied.
\end{lemma}

\begin{proof}
Restricting the field $ \xi = (\xi_t)_{t\in\R^2} =
(\xi_{t_1,t_2})_{t_1,t_2\in\R} $ to the line $ t_2=0 $ we get a process $
(\xi_{t,0})_{t\in\R} $ distributed like the process $ \xi $ of Sect.\
\ref{sect4}. On the other hand, the field $ \xi $ is distributed like the sum
$ \xi_- + \xi_+ $ of two independent fields defined by $ \xi_{a_1} =
\sum_{a_2} \xi_{a_1,a_2} $. Restricting ourselves to the line $ t_2=0 $ we see
that the process $ (\xi_{t,0})_{t\in\R} $ is distributed like the sum $
\( \xi_-(t,0) + \xi_+(t,0) \)_{t\in\R} $ of two independent processes. This
decomposition is similar to the decomposition of Sect.\ \ref{sect4} (but
different; in fact, $ \Ex \xi_-(s,0) \overline{ \xi_-(t,0) } =
\Phi((-s-t)/\sqrt2) \exp \( -0.5(s-t)^2 \) $ and $ \Ex \xi_+(s,0) \overline{
\xi_+(t,0) } = \Phi( (s+t)/\sqrt2 ) \exp \( -0.5(s-t)^2 \) $.
Still, $ \sum_k \| \Xi_- (0.25k) \| < \infty $, where $ \Xi_-(t) =
\Xi_{-,-}(t,0) + \Xi_{-,+}(t,0) $. We introduce four random fields
\[
\xi_{a_1}^{b_1} = \sum_{a_2} \xi_{a_1,a_2}^{b_1,a_2} \, ;
\]
they are independent, and each $ \xi_{a_1}^{b_1} $ is distributed like $
\xi_{a_1} $. Also,
\begin{equation}\label{F4}
\begin{aligned}
X^{0,0} &= \al \ln | \xi_-^- + \xi_+^+ | + 0.5 \al C_{\text{Euler}} \, , \\
X^{+,0} &= \al \ln | \xi_-^+ + \xi_+^+ | + 0.5 \al C_{\text{Euler}} \, , \\
X^{-,0} &= \al \ln | \xi_-^- + \xi_+^- | + 0.5 \al C_{\text{Euler}} \, .
\end{aligned}
\end{equation}
The argument of Sect.\ \ref{sect4} applies, giving
\[
\Ex \exp \bigg( \int_{-\infty}^0 | X_{t,0}^{-,0} - X_{t,0}^{0,0} | \, \D t +
\int_0^\infty | X_{t,0}^{+,0} - X_{t,0}^{0,0} | \, \D t \bigg) < \infty \, ,
\]
which is (c1); (c2) is similar.
\end{proof}

In order to verify (c3) we need a two-dimensional counterpart of Theorem
\ref{th3.1}.

\begin{proposition}\label{6.**}
Let a number $ C \in (0,\infty) $, a measurable set $ A \subset \R^2 $ and
Gaussian random fields $ \xi, \eta, \eta' $ on $ A $ be such that

(a) $ \xi, \eta, \eta' $ are independent;

(b) $ \eta $ and $ \eta' $ are identically distributed;

(c) for all $ u,v \in [0,1) $,
\[
\sum_{k,l\in\Z:(k+u,l+v)\in A} \( \Ex | \eta(k+u,l+v) |^2 \)^{1/2} \le C \, ;
\]
\indent
(d) for all $ u,v \in [0,1) $, $ n \in \{1,2,\dots\} $ and $
a_{k,l} \in \C $ ($ k,l = 0,\dots,n $),
\begin{multline*}
\Ex \bigg| \sum_{k,l\in\{0,\dots,n\}:(k+u,l+v)\in A} a_{k,l} \( \xi(k+u,l+v) +
 \eta(k+u,l+v) \) \bigg|^2 \ge \\
\ge \sum_{k,l\in\{0,\dots,n\}:(k+u,l+v)\in A} |a_{k,l}|^2 \, .
\end{multline*}
\noindent Then
\[
\Ex \exp \iint_A \ln^+ \frac{ | \xi(t_1,t_2) + \eta'(t_1,t_2) | }{ |
\xi(t_1,t_2) + \eta(t_1,t_2) | } \, \D t_1 \D t_2 \le \exp \( 2\pi (C^2+C) \)
\, .
\]
\end{proposition}

The proof is left to the reader; it is completely similar to the proof of
Theorem \ref{th3.1}.

\begin{lemma}
Condition \ref{defE3}(c3) is satisfied (for some $ \al $).
\end{lemma}

\begin{proof}
We split the given integral over $ \R^2 $ into $ 8 $ integrals according to $
\sgn t_1 $, $ \sgn t_2 $ and $ \sgn(t_1-t_2) $. These $ 8 $ integrals being
identically distributed (by symmetry), we take one of them; using the H\"older
inequality and evident symmetries,
\begin{multline*}
\Ex \exp \iint_{\R^2} | X^{0,0}_{t_1,t_2} - X^{\sgn t_1,0}_{t_1,t_2} -
 X^{0,\sgn t_2}_{t_1,t_2} + X^{\sgn t_1,\sgn t_2}_{t_1,t_2} | \, \D t_1 \D t_2
 \le \\
\le \Ex \exp 8 \iint_{t_1>t_2>0} | X^{0,0}_{t_1,t_2} - X^{+,0}_{t_1,t_2} -
 X^{0,+}_{t_1,t_2} + X^{+,+}_{t_1,t_2} | \, \D t_1 \D t_2 \le \\
\le \bigg( \Ex \exp 16 \iint_{t_1>t_2>0} | X^{0,0}_{t_1,t_2} -
 X^{+,0}_{t_1,t_2} | \, \D t_1 \D t_2 \bigg)^{1/2} \cdot \\
\cdot \bigg( \Ex \exp 16 \iint_{t_1>t_2>0} | X^{0,+}_{t_1,t_2} -
 X^{+,+}_{t_1,t_2} | \, \D t_1 \D t_2 \bigg)^{1/2} \, .
\end{multline*}
By \eqref{F4},
\[
X^{0,0}_{t_1,t_2} - X^{+,0}_{t_1,t_2} = \al \ln \frac{ | \xi_-^- (t_1,t_2) +
  \xi_+^+ (t_1,t_2) | }{ | \xi_-^+ (t_1,t_2) + \xi_+^+ (t_1,t_2) | } \, .
\]
Similarly to Sect.\ \ref{sect4},
\begin{multline*}
\Ex \exp 32 \al \iint_{t_1>t_2>0} \ln^+ \frac{ | \xi_-^- (t_1,t_2) + \xi_+^+
 (t_1,t_2) | }{ | \xi_-^+ (t_1,t_2) + \xi_+^+ (t_1,t_2) | } \, \D t_1 \D t_2 =
 \\
= \Ex \exp \iint_{t_1>t_2>0} \ln^+ \frac{ | \ti\xi (t_1,t_2) + \eta' (t_1,t_2)
 | }{ | \ti\xi (t_1,t_2) + \eta (t_1,t_2) | } \, \D t_1 \D t_2 \, ,
\end{multline*}
where $ \ti\xi, \eta $ and $ \eta' $ are defined by
\begin{align*}
\ti\xi (t_1,t_2) = 2 \xi_+^+ \Big( \frac{t_1}{\sqrt{32\al}},
 \frac{t_2}{\sqrt{32\al}} \Big) \, , \\
\eta (t_1,t_2) = 2 \xi_-^+ \Big( \frac{t_1}{\sqrt{32\al}},
 \frac{t_2}{\sqrt{32\al}} \Big) \, , \\
\eta' (t_1,t_2) = 2 \xi_-^- \Big( \frac{t_1}{\sqrt{32\al}},
 \frac{t_2}{\sqrt{32\al}} \Big) \, .
\end{align*}

Finiteness of this expectation is ensured by Prop.\ \ref{6.**} applied to
\[
A = \{ (t_1,t_2) : t_1 > t_2 > 0 \}
\]
provided that Conditions \ref{6.**}(a,b,c,d) are satisfied by $ \ti\xi $, $
\eta $ and $ \eta' $ (for some $ \al $ and $ C $).

Conditions \ref{6.**}(a,b) follow from \eqref{6.*}. Condition (c) is satisfied
since
\[
\Ex | \eta(t_1,t_2) |^2 = 2^2 \Ex \Big| \xi_-^+ \Big(
\frac{t_1}{\sqrt{32\al}}, \frac{t_2}{\sqrt{32\al}} \Big) \Big|^2 = 4 \Big\|
\Xi_- \Big( \frac{t_1}{\sqrt{32\al}}, \frac{t_2}{\sqrt{32\al}} \Big) \Big\|^2
\]
is exponentially small whenever $ t_1 \gg 1 $ (irrespective of $ t_2 $); here
$ \Xi_-(t_1,t_2) = \Xi_{-,-}(t_1,t_2) + \Xi_{-,+}(t_1,t_2) $. (The condition $
t_1 > t_2 > 0 $ is also used.)

In order to check \ref{6.**}(d) we consider (for given $ u,v \in [0,1) $)
vectors
\[
x_k = x_{k_1,k_2} = 2 \Xi \Big( \frac{k_1+u}{\sqrt{32\al}},
\frac{k_2+v}{\sqrt{32\al}} \Big)
\]
and rewrite \ref{6.**}(d) in the form
\[
\Big\| \sum_k a_k x_k \Big\|^2 \ge \sum_k |a_k|^2 \, ;
\]
here $ k $ runs over a finite subset of $ \Z^2 $. (The condition $ k_1+u >
k_2+v > 0 $ is now irrelevant.) We note that in general
\[
\Big\| \sum_k a_k x_k \Big\|^2 \ge \sum_k |a_k|^2 \Big( \|x_k\|^2 - \sum_{l\ne
k} | \ip{x_k}{x_l} | \Big)
\]
(since $ | a_k \overline{a_l} \ip{x_k}{x_l} | \le 0.5 ( |a_k|^2 + |a_l|^2 ) |
\ip{x_k}{x_l} | $), therefore the condition
\[
\| x_k \|^2 - \sum_{l\ne k} | \ip{x_k}{x_l} | \ge 1 \quad \text{for all } k
\]
is sufficient for \ref{6.**}(d).

We have
\[
| \ip{ x_{k_1,k_2} }{ x_{l_1,l_2} } | = 4 \exp \Big( -0.5 \frac{ (k_1-l_1)^2 +
(k_2-l_2)^2 }{ 32 \al } \Big) \, .
\]
Clearly, $ \sum_{l\ne k} | \ip{x_k}{x_l} | \to 0 $ as $ \al \to 0 $, uniformly
in $ k $. Choosing $ \al $ such that $ \sum_{l\ne k} | \ip{x_k}{x_l} | \le 3 $
we get $ \| x_k \|^2 - \sum_{l\ne k} | \ip{x_k}{x_l} | \ge 1 $, therefore,
\ref{6.**}(d).

\begin{sloppypar}
By Prop.\ \ref{6.**}, $ \Ex \exp 32 \iint_{t_1>t_2>0} \( X^{0,0}_{t_1,t_2} -
X^{+,0}_{t_1,t_2} \)^+ \, \D t_1 \D t_2 < \infty $. By symmetry, $ \Ex \exp
32 \iint_{t_1>t_2>0} \( X^{0,0}_{t_1,t_2} - X^{+,0}_{t_1,t_2} \)^- \, \D t_1
\D t_2 < \infty $, thus,
$ \Ex \exp 16 \iint_{t_1>t_2>0} | X^{0,0}_{t_1,t_2} -
X^{+,0}_{t_1,t_2} | \, \D t_1 \D t_2 < \infty $. Similarly, $ \Ex \exp 16
\iint_{t_1>t_2>0} | X^{0,+}_{t_1,t_2} - X^{+,+}_{t_1,t_2} | \, \D t_1 \D t_2 <
\infty $, since
\[
X^{0,+}_{t_1,t_2} - X^{+,+}_{t_1,t_2} = \al \ln \frac{ | \sum_{a_1,a_2}
\xi_{a_1,a_2}^{a_1,+} (t_1,t_2) | }{ | \sum_{a_1,a_2} \xi_{a_1,a_2}^{+,+}
(t_1,t_2) | } = \al \ln \frac{ | \zeta_-^- (t_1,t_2) +
\zeta_+^+ (t_1,t_2) | }{ | \zeta_-^+ (t_1,t_2) + \zeta_+^+ (t_1,t_2) | } \, ,
\]
where $ \zeta_{a_1}^{b_1} = \sum_{a_2} \xi_{a_1,a_2}^{b_1,+} $ are four
independent random fields, and each $ \zeta_{a_1}^{b_1} $ is distributed like
$ \xi_{a_1} $.
\end{sloppypar}
\end{proof}

\begin{remark}
Similarity between Sections \ref{sect4} and \ref{sect6} is broken near the
end; the proof of \ref{6.**}(d) is quite different from the proof of
\ref{th3.1}(d). Fourier transform works in Sect.\ \ref{sect4} but fails in
Sect.\ \ref{sect6}, since the field $ \xi $ is not stationary. This is why $
\al ( \ln |\xi| + \const ) $ is splittable for $ \al $ small enough in Sect.\
\ref{sect6}, while in Sect.\ \ref{sect4} it is splittable for $ \al=1 $ (and
in fact for all $ \al < \infty $, by the same argument). However, it does not
matter for Theorems \ref{theorem1} and \ref{theoremE}, since their conclusions
are insensitive to such coefficients.
\end{remark}

All conditions of Def.\ \ref{defE3} are verified, and we conclude.

\begin{theorem}\label{6.10}
The stationary random field $ X $ is splittable.
\end{theorem}

%% file: sect7.tex
Theorem \ref{theorem2}, Corollary \ref{corollary6} and Corollary
\ref{corollary7}, formulated in the introduction, are proved in this section.

The random entire function $ \psi(z) = \sum_{k=0}^\infty \frac{ \zeta_k z^k }{
\sqrt{k!} } $ (where $ \zeta_k $ are independent standard Gaussian) is a
centered Gaussian complex-valued random field on $ \C $ such that $ \Ex
\psi(z_1) \overline{\psi(z_2)} = \exp ( z_1 \overline z_2 ) $ for all $
z_1,z_2 \in \C $. We define another centered Gaussian complex-valued random
field $ \xi $ on $ \R^2 $ by
\[
\xi_{t_1,t_2} = \exp \( -0.5 ( |t_1|^2 + |t_2|^2 ) \) \psi (t_1 + t_2 \I)
\]
and get
\[
\Ex \xi_s \overline{ \xi_t } = \exp ( -\I s \wedge t - 0.5 |s-t|^2 ) \, ,
\]
just as in Sect.\ \ref{sect6}.

By a test function we mean a compactly supported $ C^2 $-function $ h : \R^2
\to \R $. Random variables
\[
Z(h) = \sum_{z:\psi(z)=0} h ( \Re z, \Im z )
\]
are investigated in \cite{SoTs1} and other works;
\[
\Ex Z(h) = \frac1\pi \iint_{\R^2} h(t_1,t_2) \, \D t_1 \D t_2
\]
(an immediate consequence of the Edelman-Kostlan formula), and
\[
Z(h) - \Ex Z(h) = \frac1{2\pi} \iint_{\R^2} \ln | \xi(t_1,t_2) | f(t_1,t_2) \,
\D t_1 \D t_2
\]
where $ f = \De h $, that is, $ f(t_1,t_2) = \( \frac{ \pd^2 }{ \pd t_1^2 } +
\frac{ \pd^2 }{ \pd t_2^2 } \) h(t_1,t_2) $.

Given $ r \in (0,\infty) $, we introduce
\[
h_r (t_1,t_2) = h \Big( \frac{t_1}r, \frac{t_2}r \Big) \, , \quad f_r
(t_1,t_2) = f \Big( \frac{t_1}r, \frac{t_2}r \Big)
\]
and note that
\[
\De h_r = \frac1{r^2} f_r \, .
\]
By \cite[(0.6) and Sect.~3.3]{SoTs1},
\[
\Var Z(h_r) = \frac{ \kappa }{ r^2 } \| \De h \|^2_{L_2(\R^2)} + o \Big(
\frac1{r^2} \Big) \quad \text{as } r \to \infty \, ,
\]
where $ \kappa \in (0,\infty) $ is an absolute constant.

Using $ X $ and $ \al $ of Sect.\ \ref{sect6} we have $ \ln |\xi_{t_1,t_2}| =
\frac1\al X_{t_1,t_2} + \const $, thus (taking into account that $ \iint
f_r(t_1,t_2) \, \D t_1 \D t_2 = 0 $),
\[
Z(h_r) - \Ex Z(h_r) = \frac1{2\pi\al r^2} \iint_{\R^2} f_r(t_1,t_2)
X_{t_1,t_2} \, \D t_1 \D t_2 \, .
\]
By Theorem \ref{6.10}, $ X $ is splittable. Theorem \ref{theoremE} gives $ \si
\in [0,\infty) $ (an absolute constant) such that
\[
\lim_\myatop{ r\to\infty }{ \la\log^2 r\to0 } \frac1{r^2\la^2} \ln \Ex \exp
\la \iint_{\R^2} f_r(t_1,t_2) X_{t_1,t_2} \, \D t_1 \D t_2 = \frac{ \si^2 }2
\| f \|^2_{L_2(\R^2)} \, ,
\]
that is,
\[
\lim_\myatop{ r\to\infty }{ \la\log^2 r\to0 } \frac1{r^2\la^2} \ln \Ex \exp
2\pi\al r^2 \la \( Z(h_r) - \Ex Z(h_r) \) = \frac{ \si^2 }2 \| f \|^2 \, ,
\]
or equivalently,
\[
\lim_\myatop{ r\to\infty }{ \la\log^2 r\to0 } \frac1{r^2\la^2} \ln \Ex \exp
r^2 \la \( Z(h_r) - \Ex Z(h_r) \) = \frac{ \si^2 }{2 (2\pi\al)^2} \| f \|^2 \,
.
\]
It follows that
\begin{multline*}
\frac{ \si^2 }{2 (2\pi\al)^2} \| f \|^2 = \lim_{r\to\infty} \frac1{r^2}
 \lim_{\la\to0} \frac1{\la^2} \ln \Ex \exp r^2 \la \( Z(h_r) - \Ex Z(h_r) \) =
 \\
= \lim_{r\to\infty} \frac1{r^2} \cdot \frac12 r^4 \Var Z(h_r) = \frac\kappa2
 \|f\|^2 \, ,
\end{multline*}
that is, $ \si = 2\pi\al \sqrt\kappa $. We get
\[
\lim_\myatop{ r\to\infty }{ \la\log^2 r\to0 } \frac1{r^2\la^2} \ln \Ex \exp
r^2 \la \( Z(h_r) - \Ex Z(h_r) \) = \frac\kappa2 \| f \|^2 \, ,
\]
which proves Theorem \ref{theorem2}.

By Corollary \ref{5.3},
\[
\lim_\myatop{ r\to\infty, c\to\infty }{ (c\log^2 r)/r \to0 } \frac1{c^2} \ln
\PR{ \iint f_r (t_1,t_2) X_{t_1,t_2} \, \D t_1 \D t_2 \ge c\si \| f \| r } =
-\frac12 \, ,
\]
that is,
\[
\lim_\myatop{ r\to\infty, c\to\infty }{ (c\log^2 r)/r \to0 } \frac1{c^2} \ln
\PR{ 2\pi\al r^2 \( Z(h_r) - \Ex Z(h_r) \) \ge c\si \| f \| r } = -\frac12 \,
.
\]
Taking into account that $ \si = 2\pi\al \sqrt\kappa $ we get
\[
\lim_\myatop{ r\to\infty, c\to\infty }{ (c\log^2 r)/r \to0 } \frac1{c^2} \ln
\PR{ Z(h_r) - \Ex Z(h_r) \ge \frac{ c\sqrt\kappa \| f \| }{r} } = -\frac12 \,
,
\]
which proves Corollary \ref{corollary6}.

By Corollary \ref{5.4}, the distribution of $ r^{-1} \iint f_r (t_1,t_2)
X_{t_1,t_2} \, \D t_1 \D t_2 $ converges (as $ r \to \infty $) to the normal
distribution $ N(0,\si^2 \| f \|^2) $. That is,
\[
2\pi\al r \( Z(h_r) - \Ex Z(h_r) \) \to N(0,\si^2 \| f \|^2) \quad \text{in
distribution.}
\]
and therefore
\[
r \( Z(h_r) - \Ex Z(h_r) \) \to N \bigg( 0, \Big( \frac{\si \| f \|}{2\pi\al}
\Big)^2 \bigg) = N(0,\kappa \| f \|^2 ) \quad \text{in distribution,}
\]
which proves Corollary \ref{corollary7}.

%% file: main.bbl
\begin{thebibliography}{8.}

{\raggedright
\bibitem{Bi} P. Billinglsley (1995):
\emph{Probability and measure} (third edition),
Wiley.

\bibitem{Bo} V.I. Bogachev (1998):
\emph{Gaussian measures,}
AMS.

\bibitem{DMPU} J. Dedecker, F. Merlevede, M. Peligrad, S. Utev (2007):
\emph{Moderate deviations for stationary sequences of bounded random
 variables,}
\texttt{arXiv:0711.3924}.

\bibitem{DGW} H. Djellout, A. Guillin, L. Wu (2006):
\emph{Moderate deviations of empirical periodogram and non-linear functionals
 of moving average process,}
Ann. Inst. H. Poincar\'e Probab. Statist. \textbf{42}:4, 393--416.

\bibitem{Ell} R.S. Ellis (2006):
\emph{The theory of large deviations and applications to statistical
 mechanics,}
\href{http://www.math.umass.edu/~rsellis/pdf-files/Dresden-lectures.pdf}%
{http://www.math.umass.edu/$\sim$rsellis/pdf-files/Dresden-lectures.pdf}

\bibitem{KS} J.R. Klauder, E.C.G. Sudarshan (1968):
\emph{Fundamentals of quantum optics,}
W.A. Benjamin, Inc.

\bibitem{NSV} F. Nazarov, M. Sodin, A. Volberg (2007):
\emph{The
Jancovici\nobreakdash-\hspace{0pt}Lebowitz\nobreakdash-\hspace{0pt}Manificat
law for large fluctuations of random complex zeroes,}
\texttt{arXiv:0707.3863}.

\bibitem{SoTs1} M. Sodin, B. Tsirelson (2004):
\emph{Random complex zeroes, I. Asymptotic normality,}
Israel Journal of Mathematics \textbf{144}, 125--149.
Also, \texttt{arXiv:math.CV/0210090}.

}
\end{thebibliography}
